\documentclass[
submission]{pmsart}[2006/01/13]

\usepackage{color}

\input cyracc.def

    \def\one{{\mathbb I}}
    \def\BbbE{\mathbb E}
    \def\BbbN{\mathbb N}

    \def\cB{\mathcal{B}}
    \def\cF{\mathcal F}
    \def\cS{\mathcal S}
    \def\calR{\mathcal R}
    \def\cT{\mathcal T}

    \def\bE{\mathbf E}
    \def\bH{\mathbf H}
    \def\bP{\mathbf P}
    \def\bT{\mathbf T}
    \def\bQ{\mathbf Q}

        \def\vecx{\vec{x}}
        \def\vecX{\vec{X}}
        \def\vecalpha{\vec{\alpha}}
        \def\vecPi{\overrightarrow{\Pi}}
        
        \def\plus{+}
\def\esssup{\operatornamewithlimits{ess\,sup}}

\def\rightmaltese{\protect\vspace*{-4ex}
\begin{flushright}\(\maltese\)\end{flushright}}
\newenvironment{pf}{{\sc Proof.}\hspace{0.5cm}}{\rightmaltese}
\renewcommand{\baselinestretch}{1.4}

\volume{0}%
\fasc{0}%
\years{20xx}%
\pages{000--000}%
\firstpage{1} %
\lastrevision{2nd of January 2010}
\title{On a random number of disorders}
\titlethanks{}
\dedicated{}
\ShortTitle{Random number of disorders} 
\ShortAuthors{K. Szajowski}
\NumberOfAuthors{1}

\received{23rd of November 2008}
\FirstAuthor{Krzysztof Szajowski}
\FirstAuthorCity{Wroc\l{}aw}
\FirstAuthorAffiliation{Institute of Mathematics and Computer Science}
\FirstAuthorAddress{Wroc\l{}aw University of Technology\\
                    Wybrze\.{z}e Wyspia\'{n}skiego 27\\
                    50-370 Wroc\l{}aw, Poland}
\FirstAuthorEmail{Krzysztof.Szajowski@pwr.wroc.pl}                              
\FirstAuthorThanks{Affiliation: Institute of Mathematics, Polish Academy of Science, \'Sniadeckich~8, 00-956~Warszawa, Poland}                             

\SecondAuthor{} 
\SecondAuthorCity{} %
\SecondAuthorAffiliation{}%
\SecondAuthorAddress{}%
\SecondAuthorEmail{}%
\SecondAuthorThanks{}%

\ThirdAuthor{}
\ThirdAuthorCity{}%
\ThirdAuthorAffiliation{}%
\ThirdAuthorAddress{}%
\ThirdAuthorEmail{}
\ThirdAuthorThanks{}%

\FourthAuthor{}
\FourthAuthorCity{}%
\FourthAuthorAffiliation{}%
\FourthAuthorAddress{}%
\FourthAuthorEmail{}%
\FourthAuthorThanks{}%

\FifthAuthor{} 
\FifthAuthorCity{} %
\FifthAuthorAffiliation{}%
\FifthAuthorAddress{} %
\FifthAuthorEmail{}%
\FifthAuthorThanks{}%


\MSCcodes{Primary: 60G40 60K99; Secondary: 90D60.}               

\begin{document}
\bigskip
\maketitle
\medskip

\begin{abstract}
\baselineskip15pt
We register a random sequence which has the following properties: it has three segments being the homogeneous Markov processes. Each segment has his own one step transition probability law and the length of the segment is unknown and random. It means that at two random moments $\theta_1$, $\theta_2$, where $0\leq \theta_1 \leq \theta_2$, the source of observations is changed and the first observation in new segment is chosen according to new transition probability starting from the last state of the previous segment. In effect the number of homogeneous segments is random. The transition probabilities of each process are known and \emph{a priori} distribution of the disorder moments is given. The former research on such problem has been devoted to various questions concerning the distribution changes. The random number of distributional segments creates new problems in solutions with relation to analysis of the model with deterministic number of segments. Two cases are presented in details. In the first one the objectives is to stop on or between the disorder moments while in the second one our objective is to find the strategy which immediately detects the distribution changes. Both problems are reformulated to optimal stopping of the observed sequences. The detailed analysis of the problem is presented to show the form of optimal decision function.
\end{abstract}
\renewcommand{\baselinestretch}{1.4}
\section{Introduction}
\label{wstep}
Suppose that the process $X=\{X_n,n\in\BbbN\}$, $\BbbN=\{0,1,2,\ldots\}$, is observed sequentially. The process is obtained from three Markov processes by switching between them at two random moments of time, $\theta_1$ and $\theta_2$. Our objective is to detect these moments based on observation of $X$.

Such model of data appears in many practical problems of the quality control~(see Brodsky and Darkhovsky~\cite{brodar93:nonparametr}, Shewhart~\cite{she31:quality} and in the collection of the papers \cite{basben86:abrupt}), traffic anomalies in networks (in papers by Dube and Mazumdar~\cite{dubmaz01:quickest}, Tartakovsky et al.~\cite{tarroz06:intrusions}), epidemiology models (see Baron~\cite{bar04:epidemio}). In management of manufacture it happens that the plants which produce some details changes their parameters. It makes that the details change their quality. Production can be divided into three sorts. Assuming that at the beginning of production process the quality is highest, from some moment $\theta_1$ the products should be classified to lower sort and beginning with the moment $\theta_2$ the details should be categorized as having the lowest quality. The aim is to recognize the moments of these changes.

Shiryaev~\cite{shi61:detection,shi78:optimal} solved the disorder problem of the independent random variables with one disorder where the mean distance between disorder time and the moment of its detection was minimized. The probability maximizing approach to the problem was used by Bojdecki~\cite{boj79:disorder} and the stopping time which is in a given neighborhood of the moment of disorder with maximal probability was found. The disorders in more complicated dependence structures of switched sequences are subject of investigation by Pelkowitz~\cite{pel87:discrete,pel87:Markov}, Yakir~\cite{yak94:finite}, Mustakides~\cite{mou98:abrupt}, Lai~\cite{lai95:changepoint,lai98:quick}, Fuh~\cite{fuh04:hidden}, Tartakovsky and Veeravalli~\cite{tarvee08:quickest}. The probability maximizing approach to such problems with two disorders was considered by Yoshida~\cite{yos83:complicated}, Szajowski~\cite{sza92:detection,sza96:twodis} and Sarnowski and Szajowski~\cite{sarsza08:disorder}. Yoshida~\cite{yos83:complicated} investigated the problem of optimal stopping the observation of the process $X$ so as to maximize the probability that the distance between the moments of disorder $\theta_i$ and their estimates, the stopping times $\tau_i$, $i=1,2$, will not exceed given numbers (for each disorder independently). This question has been reformulated by Szajowski~\cite{sza96:twodis} to the simultaneous detection of both disorders under requirement that the performance of procedure is globally measured for both detections and it has been extended to the case with unknown distribution between disorders by Sarnowski and Szajowski~\cite{sarsza08:disorder} (see also papers by Bojdecki and Hosza~\cite{bojhos84:problem} for related approach with switching sequences of independent random variables). The method of solution is based on a~transformation of the model to the double optimal stopping problem for markovian function of some statistics (see Haggstrom~\cite{hag67:more}, Nikolaev~\cite{nik79:obob}). The strategy which stops the process between the first and the second disorder with maximal probability has been constructed by Szajowski~\cite{sza92:detection}. The considerations are inspired by the problem regarding how we can protect ourselves against a second fault in a technological system after the occurrence of an initial fault or by the problem of detection the beginning and the end of an epidemic.

The paper is devoted to a generalization of the double disorder problem considered both in~\cite{sza92:detection} and
\cite{sza96:twodis} in which immediate switch from the first preliminary distribution to the third one is possible (i.e. it is possible that the random variables $\theta_1$ and $\theta_2$ are equal with a positive probability). It is also possible that we observe the homogeneous data without disorder when both disorder moments are equal to $0$. The extension leads to serious difficulties in the construction of an equivalent double optimal stopping model. The formulation of the problem can be found in Section \ref{sformProblem}. The main results are subject of Sections \ref{SolDetection1} (see Theorem~\ref{twierdz1}) and
\ref{SolDisorders1}.

\section{Formulation of detection problems}\label{sformProblem}
\vspace{-.5cm}
Let $(X_n)_{ n \in \BbbN}$ be an observable sequence of random variables defined on the space $(\Omega,\mathcal{F}, \bP)$ with values in $(\BbbE, \mathcal{B})$, where $\BbbE$ is a Borel subset of $\mathbf{R}$. On $(\BbbE, \mathcal{B})$ there is $\sigma$-additive measure $\mu$. On the same probability space there are defined random variables $\theta_1$, $\theta_2$ with values in $\BbbN$ and the following distributions:
\begin{eqnarray}
\label{rozkladyTeta}
\bP(\theta_1 = j) &=& \one_{\{j=0\}}(j)\pi+\one_{\{j>0\}}(j)\bar{\pi} p_1^{j-1}q_1,\\
\label{rokladWarTeta2}
\bP(\theta_2 = k \mid \theta_1=j) &=&\one_{\{k=j\}}(k)\rho+\one_{\{k>j\}}(k)\bar{\rho} p_2^{k-j-1}q_2
\end{eqnarray}
where $j=0,1,2,...$, $k=j,j+1,j+2,...$, $\bar\pi=1-\pi$, $\bar\rho=1-\rho$.
Additionally we consider Markov processes $(X_n^{i}, \mathcal{G}_n^{i}, \bP_x^{i})$ on $(\Omega,\mathcal{F}, \bP)$,
$i=0,1,2$, where $\sigma$-fields $\mathcal{G}_n^{i}$ are the smallest $\sigma$-fields for which
 ${(X^{i}_n)}_{n=0}^\infty$, $i=0,1,2$, are adapted, respectively. Let us define process $(X_n)_{ n \in \BbbN}$ in the following way:
\begin{eqnarray}
\label{procesyX}
    X_n = X^{0}_n  \one_{\{\theta_1>n\}} + X^{1}_n\one_{\{X^{1}_{\theta_1-1}=X^{0}_{\theta_1-1},\theta_1 \leq n < \theta_2\}} + X^{2}_n\one_{\{X^{2}_{\theta_2-1}=X^{1}_{\theta_2-1},\theta_2 \leq n\}}.
\end{eqnarray}
We make inference on $\theta_1$ and $\theta_2$ from the observable sequence ($X_n$, $n \in \BbbN$) only. %
It should be emphasized that the sequence ($X_n$, $n \in \BbbN$) is not markovian under admitted assumption as it has been mentioned in \cite{sza92:detection}, \cite{yak94:finite} and \cite{dubmaz01:quickest}. However, the sequence satisfies the Markov property given $\theta_1$ and $\theta_2$ (see Szajowski~\cite{sza96:twodis} and Moustakides~\cite{mou98:abrupt}). Thus for further consideration we define filtration $\{\mathcal{F}_n\}_{n \in \BbbN}$, where $\mathcal{F}_n = \sigma(X_0,X_1,...,X_n)$, related to real observation. Variables $\theta_1$, $\theta_2$ are not stopping times with respect to $\mathcal{F}_n$ and $\sigma$-fields $\mathcal{G}_n^{\bullet}$. Moreover, we have knowledge about the distribution of $(\theta_1,\theta_2)$ independent of any observation of the sequence $(X_n)_{n \in \BbbN}$. This distribution, called the \emph{a priori distribution} of $(\theta_1,\theta_2)$ is given by (\ref{rozkladyTeta}) and (\ref{rokladWarTeta2}). 

It is assumed that the measures ${\bP^i_x}(\cdot)$ on $\cF$, $i=0,1,2$, have following representation. For any $B\in\cB$ we have 
\[
\bP^i_x(\omega:X_1^i\in B)=\bP(X_1^i\in B|X_0^i=x)=\int_B f_x^i(y)\mu(dy)=\int_B \mu_x^i(dy)=\mu^i_x(B),
\] 
where the functions $f_x^{i}(\cdot)$ are different and $f_x^{i}(y)/f_x^{(i+1)\text{mod} 3}(y) < \infty$ for $i=0,1,2$ and all $x,y \in \BbbE$. We assume that the measures $\mu_x^{i}$, $x \in \BbbE$ are known in advance.

For any $D_n=\{\omega:X_i\in B_i,\; i=1,\ldots,n\}$, where $B_i\in \cB$, and any $x\in\BbbE$ define
\[
\bP_x(D_n)=\bP(D_n|X_0=x)=\int_{\times_{i=1}^n B_i}S_n(x,\vec{y}_n) \mu(d\vec{y}_n)=\int_{\times_{i=1}^nB_i} \mu_x(d\vec{y}_n)=\mu_x(\times_{i=1}^nB_i), 
\] 
where the sequence of functions $S_n:\times_{i=1}^n\BbbE\rightarrow\Re$ is given by (\ref{disXuncond}) in Appendix.

The presented model has the following heuristic justification: two disorders take place in the observed sequence $(X_n)$. They affect distributions by changing their parameters. The disorders occur at two random times $\theta_1$ and $\theta_2$, $\theta_1\leq \theta_2$. They split the sequence of observations into segments, at most three ones. The first segment is described by $(X_n^0)$, the second one - for $\theta_1 \leq n < \theta_2$ - by $(X_n^1)$. The third is given by $(X_n^2)$ and is observed when $n\geq \theta_2$. When the first disorder takes the place there is a "switch" from the initial distribution to the distribution with the conditional density $f_x^i$ with respect of the measure $\mu$, where $i=1$ or $i=2$, when $\theta_1<\theta_2$ or $\theta_1=\theta_2$, respectively. Next, if $\theta_1<\theta_2$, at the random time $\theta_2$ the distribution of observations becomes $\mu_x^2$. We assume that the variables $\theta_1, \theta_2$ are unobservable directly.

Let $\cS$ denote the set of all stopping times with respect to the filtration $(\cF_n)$, $n=0,1,\ldots$ and
$\cT=\{(\tau,\sigma): \tau\leq \sigma,\ \tau,\sigma\in\cS\}$. Two problems with three distributional segments are recalled to investigate them under weaker assumption that there are at most three homogeneous segments.
\subsection{Detection of change}
Our aim is to stop the observed sequence between the two disorders.This can be interpreted as a strategy for protecting against a second failure when the first has already happened. The mathematical model of this is to control the probability $\bP_x(\tau < \infty, \theta_1 \leq \tau < \theta_2)$ by choosing the stopping time $\tau^{*}\in \cS$ for which 
\vspace{-1.67ex}
\begin{equation}
\label{problem}
\bP_x(\theta_1 \leq \tau^{*} < \theta_2)
= \sup_{\tau \in \mathcal{T}}\bP_x(\tau < \infty, \theta_1 \leq \tau < \theta_2).
\end{equation}
\subsection{Disorders detection}
Our aim is to indicate the moments of switching with given precision $d_1,d_2$ (Problem $\mbox{D}_{d_1d_2}$). We want to determine a pair of stopping times $(\tau^*,\sigma^*)\in\cT$
such that for every $x\in\BbbE$
\small
\begin{equation}
\label{equ3}
\bP_x(|\tau^*-\theta_1|\leq d_1,|\sigma^*-\theta_2|\leq d_2)=
\sup_{\stackrel{(\tau,\sigma)\in\cT}{0\leq\tau\leq\sigma<\infty}} 
\bP_x(|\tau-\theta_1|\leq d_1,|\sigma-\theta_2|\leq d_2).
\end{equation}
\normalsize
The problem has been considered in \cite{sza96:twodis} under natural simplification that there are three segments of data (\emph{i.e.} there is $0<\theta_1<\theta_2$). In the section \ref{SolDisorders1} the problem $\mbox{D}_{00}$ is analyzed.

\section{On some \emph{a posteriori} processes}
\label{PosterioriProc}

The formulated problems are translated to the optimal stopping problems for some Markov processes. The important part
of the reformulation process is choice of the \emph{statistics} describing knowledge of the decision maker.
The \emph{a posteriori} probabilities of some events play the crucial role. Let us define the following
\emph{a posteriori} processes (cf. \cite{yos83:complicated}, \cite{sza92:detection}).
\begin{eqnarray}
\label{pi1x}
\Pi^i_n&=&\bP_x(\theta_i\leq n|\cF_n),\\
\label{pi12nmx}
\Pi^{12}_{n}&=&\bP_x(\theta_1=\theta_2>n|\cF_n)=P_x(\theta_1=\theta_2>n|\cF_{mn}),\\
\label{pinmx}
\Pi_{mn}&=&\bP_x(\theta_1=m,\theta_2>n|\cF_{mn}),
\end{eqnarray}
where $\cF_{m\; n}=\cF_n$ for $m,n=1,2,\ldots$, $m<n$, $i=1,2$. For recursive representation of (\ref{pi1x})--(\ref{pinmx}) we need the following functions:
\begin{eqnarray*}
\Pi^1(x,y,\alpha,\beta,\gamma)&=&1-\frac{p_1(1-\alpha) f^{0}_x(y)}{\bH(x,y,\alpha,\beta,\gamma)}\\
\Pi^2(x,y,\alpha,\beta,\gamma)&=&\frac{(q_2\alpha + p_2\beta +
q_1\gamma) f^{2}_x(y)}{\bH(x,y,\alpha,\beta,\gamma)}\\
\Pi^{12}(x,y,\alpha,\beta,\gamma)&=&\frac{p_1\gamma f^0_x(y)}{\bH(x,y,\alpha,\beta,\gamma)}\\
\Pi(x,y,\alpha,\beta,\gamma,\delta)&=&\frac{p_2\delta f^{1}_x(y)}{\bH(x,y,\alpha,\beta,\gamma)}
\end{eqnarray*}
where $\bH(x,y,\alpha,\beta,\gamma)= (1-\alpha) p_1 f^{0}_x(y)+ [p_2(\alpha -\beta) + q_1(1- \alpha-\gamma)]f^{1}_x(y)+ [q_2\alpha+p_2\beta + q_1\gamma]f^{2}_x(y)$. In the sequel we adopt the following denotations
\begin{eqnarray}
\label{vecalpha} \vecalpha&=&(\alpha,\beta,\gamma)\\
\label{vecPi}     \vecPi_n&=&(\Pi^1_n,\Pi^2_n,\Pi^{12}_n).
\end{eqnarray}
The basic formulae used in the transformation of the disorder problems to the stopping problems are given in the following
\begin{lemma}
\label{reqform}
For each $x\in\BbbE$ the following formulae,  for $m,n=1,2,\ldots$, $m<n$, hold:
\begin{eqnarray}
\label{eqpi1}
 \Pi^1_{n+1}&=&\Pi^1(X_n,X_{n+1},\Pi^1_n,\Pi^2_n,\Pi^{12}_n)\\
\label{eqpi2}
\Pi^2_{n+1}&=&\Pi^2(X_n,X_{n+1},\Pi^1_n,\Pi^2_n,\Pi^{12}_n)\\
 \label{eqpi12}
 \Pi^{12}_{n+1}&=&\Pi^{12}(X_n,X_{n+1},\Pi^1_n,\Pi^2_n,\Pi^{12}_n)\\
 \label{eqpi}
 \Pi_{m\,n+1}&=&\Pi(X_n,X_{n+1},\Pi^1_n,\Pi^2_n,\Pi^{12}_n,\Pi_{m\,n})
\end{eqnarray}
with boundary condition $\Pi^1_0=\pi$, $\Pi^2_0(x)=\pi\rho$, $\Pi^{12}_0(x)=\bar{\pi}\rho$, and 
$\Pi_{m\,m}=(1-\rho){\frac{q_1 f^1_{X_{m-1}}(X_m)}{p_1 f^0_{X_{m-1}}(X_m)}}(1-\Pi^1_m)$.
\end{lemma}

\begin{pf}
The cases (\ref{eqpi1}), (\ref{eqpi2}) and (\ref{eqpi}), when $0<\theta_1<\theta_2$, have been proved in \cite{yos83:complicated} and \cite{sza92:detection}. Let us assume $0\leq \theta_1\leq \theta_2$ and suppose that
$B_i \in {\cB}$, $1\leq i \leq n+1$.  Let us assume that $X_0=x$ and denote $D_n=\{\omega:X_i(\omega)\in B_i, 1\leq i\leq n\}$.
\begin{description}
\item[Ad. \eqref{eqpi1}] For $A_i=\{\omega:X_i\in B_i\}\in\cF_i$, $1\leq i\leq n+1$ and $D_{n+1}\in\cF_{n+1}$ we have by properties of $S_n(\vecx_{n})$ where $\vecx_n=(x_0,\ldots,x_n)$ (see Lemma~\ref{recSn})
\begin{eqnarray*}
\int_{D_{n+1}}\bP_x(\theta_1>n+1|\cF_{n+1})d\bP_x&=&\int_{D_{n+1}}\one_{\{\theta_1>n+1\}}d\bP_x\\
&\hspace{-24em}=&\hspace{-12em} \int_{\times_{i=1}^{n+1}B_i} \frac{(f^{n<\theta_1<\theta_2}_x(\vecx_{1,n})+f^{n<\theta_1=\theta_2}_x(\vecx_{1,n}))}{S_n(\vecx_{n})} \frac{p_1f^0_{x_n}(x_{n+1})}{\bH(x_n,x_{n+1},\vecPi_n(\vecx_{n}))}\mu_x(d\vecx_{1,n+1})\\
&\hspace{-24em}=&\hspace{-12em}\int_{D_{n+1}}(1-\Pi^1_n)\frac{p_1f^0_{X_n}(X_{n+1})}{\bH(X_n,X_{n+1},\vecPi_n)}d\bP_x.
\end{eqnarray*}
\setlength\arraycolsep{2pt}
Thus, taking into account (\ref{pi1x}) we have $\Pi_{n+1}^{1} = 1 - P_{x}\left(\theta_1 > n + 1 \mid\cF_{n+1} \right)= 1- (1-\Pi_{n}^{1})p_1f^0_{X_n}(X_{n+1})\bH^{-1}(X_n,X_{n+1},\vecPi_n)$. This proves the form of the formula (\ref{eqpi1}).

\item[Ad. \eqref{eqpi2}] Under the same denotations like in the proof of (\ref{eqpi1}) we have using denotation from Section~\ref{distDISsample} and the results of Lemma~\ref{multidistform}
\begin{eqnarray*}
  \int_{D_{n+1}}\bP_{x}(\theta_2 \leq n + 1 \mid \cF_{n+1})d\bP_x &=&\int_{D_{n+1}}\one_{\{\theta_2\leq n+1\}}d\bP_x\\
&\hspace{-24em}\quad\stackrel{(\ref{theta1lesstheta2lessn})}{=}&\hspace{-12em} \int_{\times_{i=1}^{n+1}B_i} \frac{f^{\theta_1\leq\theta_2\leq n+1}_x(\vecx_{1,n+1})}{S_n(\vecx_{n}) \bH(x_n,x_{n+1},\vecPi_n(\vecx_{n}))}\mu_x(d\vecx_{1,n+1})\\
&\hspace{-24em}=&\hspace{-12em} \int_{\times_{i=1}^{n+1}B_i} \frac{[q_2\Pi^1_n(\vecx_{0,n})+p_2\Pi^2_n(\vecx_{0,n})+q_1\Pi^{12}_n(\vecx_{0,n})]f^2_{x_n}(x_{n+1})}{ \bH(x_n,x_{n+1},\vecPi_n(\vecx_{n}))}\mu_x(d\vecx_{1,n+1})\\
&\hspace{-24em}=&\hspace{-12em} \int_{D_{n+1}} \frac{[q_2\Pi^1_n+p_2\Pi^2_n+q_1\Pi^{12}_n]f^2_{X_n}(X_{n+1})}{ \bH(X_n,X_{n+1},\vecPi_n)}d\bP_x.
\end{eqnarray*}
Thus we get:
\begin{eqnarray*}
  \Pi_{n+1}^2 &=& \bP_{x}(\theta_2 \leq n + 1 \mid \cF_{n+1})\\
    &=& \left[ (\Pi_n^{1} - \Pi_n^2)q_2 + \Pi_n^2+q_1\Pi_n^{12}\right]f^2_{X_n}(X_{n+1})   \bH^{-1}(X_n,X_{n+1},\vecPi_n) \nonumber
\end{eqnarray*}
which leads to the formula (\ref{eqpi2}).

\item[Ad. \eqref{eqpi12}] By (\ref{pi12nmx}) and the results of Lemma~\ref{multidistform}
\begin{eqnarray*}
  \int_{D_{n+1}}\bP_{x}(\theta_2=\theta_1 > n + 1 \mid \cF_{n+1})d\bP_x &=&\int_{D_{n+1}}\one_{\{\theta_2=\theta_1\geq n+1\}}d\bP_x\\
&\hspace{-24em}=&\hspace{-12em} \int_{\times_{i=1}^{n+1}B_i} \frac{f^{\theta_1=\theta_2> n}_x(\vecx_{1,n+1})}{S_n(\vecx_{n}) \bH(x_n,x_{n+1},\vecPi_n(\vecx_{n}))}\mu_x(d\vecx_{1,n+1})\\
&\hspace{-24em}=&\hspace{-12em} \int_{\times_{i=1}^{n+1}B_i} \frac{\Pi^{12}_n(\vecx_{n})p_1f^0_{x_n}(x_{n+1})}{ \bH(x_n,x_{n+1},\vecPi_n(\vecx_{n}))}\mu_x(d\vecx_{1,n+1})\\
&\hspace{-24em}=&\hspace{-12em} \int_{D_{n+1}} \frac{\Pi^{12}_n p_1 f^0_{X_n}(X_{n+1})}{\bH(X_n,X_{n+1},\vecPi_n)}d\bP_x,
\end{eqnarray*}
which leads to:
\begin{eqnarray*}
  \Pi_{n+1}^{12} &=& p_1\Pi_n^{12}f^0_{X_n}(X_{n+1})   \bH^{-1}(X_n,X_{n+1},\vecPi_n) \nonumber
\end{eqnarray*}
and it proves the formula (\ref{eqpi12}).

\item[Ad. \eqref{eqpi}] Similarly, by the definition (\ref{pinmx}) and the results of Lemma~\ref{multidistform} we get 
\begin{eqnarray*}
  \int_{D_{n+1}}\bP_{x}(\theta_1=m,\theta_2 > n + 1 \mid \cF_{n+1})d\bP_x &=&\int_{D_{n+1}}\one_{\{\theta_1=m,\theta_2> n+1\}}d\bP_x\\
&\hspace{-28em}=&\hspace{-14em} \int_{\times_{i=1}^{n+1}B_i} \frac{\bar{\pi}\bar{\rho}p_1^{m-1}q_1p_2^{n+1}\prod_{s=1}^{m-1}f^0_{x_{s-1}}(x_s)\prod_{k=m}^nf^1_{x_{k-1}}(x_{k})f^1_{x_n}(x_{n+1})}{S_n(x_{0,n}) \bH(x_n,x_{n+1},\vecPi_n(\vecx_{n}))}\mu_x(d\vecx_{1,n+1})\\
&\hspace{-28em}=&\hspace{-14em} \int_{\times_{i=1}^{n+1}B_i} \frac{\Pi_{m\; n}(\vecx_{n})p_2f^1_{x_n}(x_{n+1})}{ \bH(x_n,x_{n+1},\vecPi_n(\vecx_{n}))}\mu_x(d\vecx_{1,n+1})= \int_{D_{n+1}} \frac{\Pi_{m\; n} p_2 f^1_{X_n}(X_{n+1})}{\bH(X_n,X_{n+1},\vecPi_n)}d\bP_x.
\end{eqnarray*}
It leads to relation
\begin{eqnarray*}
  \Pi_{m\;n+1} &=& p_2\Pi_{m\; n}f^1_{X_n}(X_{n+1})   \bH^{-1}(X_n,X_{n+1},\vecPi_n) \nonumber
\end{eqnarray*}
and it proves the formula (\ref{eqpi}).

\end{description}
Further details concerning recursive formula for conditional probabilities can be found in Remark~\ref{reccondprobab} in Appendix. 
\end{pf}

\begin{remark}
Let us assume that the considered Markov processes have the finite state space and $\vecx_n=(x_0,x_1,\ldots,x_n)$, $x_0=x$ are given. In this case the formula (\ref{eqpi}) follows from the Bayes formula:
\[
\bP_x(\theta_1=j,\theta_2=k|\vec{X}_n=\vecx_n) =
\left\{
\begin{array}{ll}
p_{jk}^\theta\prod_{s=1}^nf^0_{x_{s-1}}(x_s)(S_n(\vecx_n))^{-1}&\mbox{ if $j>n$},\\
p_{jk}^\theta\prod_{s=1}^{j-1}f^0_{x_{s-1}}(x_s) &\mbox{}\\
\quad\times \prod_{t=j}^{n}f^1_{x_{t-1}}(x_t) (S_n(\vecx_n))^{-1} &\mbox{ if $j\leq n<k$},\\
p_{jk}^\theta\prod_{s=1}^nf^0_{x_{s-1}}(x_s)\prod_{t=j}^{k-1}f^1_{x_{t-1}}(x_{t})&\mbox{}\\
\quad\times\prod_{u=k}^{n}f^2_{x_{u-1}}(x_u)(S_n(\vecx_n))^{-1} &\mbox{ if $k\leq n$},
\end{array}
\right.
\]
where $p_{jk}^\theta=\bP(\theta_1=j,\theta_2=k)$ and $S_n(\cdot)$ is given by (\ref{disXuncond}).
\end{remark}

\begin{lemma}\label{lematTechniczny}
For each $x\in\BbbE$ and each Borel function $u:\BbbE\longrightarrow \Re$ the following equations are fulfilled:
\small
\begin{eqnarray}
\label{wartOczPi1}
&&\bE_{x}\left(u(X_{n+1})(1 - \Pi_{n+1}^{1})\mid \mathcal{F}_n\right) = (1-\Pi^{1}_n)p_1\int_{\BbbE}u(y)f^0_{X_n}(y)\mu(dy),\\
\label{wartOczPi1Pi2}
&&\bE_{x}\left(u(X_{n+1})(\Pi_{n+1}^{1} - \Pi_{n+1}^2)\mid \cF_n\right)   \\
    &&\quad= \left[q_1(1-\Pi^{1}_n-\Pi^{12}_n)+p_2(\Pi^{1}_n-\Pi^{2}_n) \right]
    \int_{\BbbE}u(x)f^1_{X_n}(y)\mu(dy),\nonumber\\
\label{wartOczPi2}
&&\bE_{x}\left(u(X_{n+1})\Pi_{n+1}^{2})\mid \cF_n\right) = \left[ q_2\Pi^1_n+p_2\Pi_n^2+q_1\Pi^{12}_n \right]\!\int_{\BbbE}\!u(y)f^2_{X_n}(y)\mu(dy),\\
\label{wartOczPi12}
\mbox{$\;$}&&\bE_{x}\left(u(X_{n+1})\Pi^{12}_{n+1})\mid \cF_n\right) =p_1\Pi^{12}_n\int_{\BbbE}u(y)f^0_{X_n}(y)\mu(dy)
\end{eqnarray}
\normalsize
\begin{equation}
\label{exp1}
\bE_x(u(X_{n+1})|\cF_n) = \int_{\BbbE} u(y)\bH(X_n,y,\vecPi_n)\mu(dy).
\end{equation}
\end{lemma}
\begin{pf}
The relations \eqref{wartOczPi1}-\eqref{wartOczPi12} are consequence of suitable division of $\Omega$ defined by
$(\theta_1,\theta_2)$ and properties established in Lemma~\ref{multidistform}. Let us prove the equation (\ref{wartOczPi2}). To this end let us define $\sigma$-field $\mathcal{\widetilde{F}}_{n} = \sigma(\theta_1,\theta_2,X_0,...,X_n)$.
Notice that $\cF_{n} \subset \mathcal{\widetilde{F}}_{n}$. We have:
\begin{eqnarray*}
\bE_{x}(u(X_{n+1})\Pi_{n+1}^{2}\mid \mathcal{F}_{n}) &=& 
\bE_{x}(u(X_{n+1})\bE_x(\one_{\{\theta_2\leq n+1\}}\mid\cF_{n+1})\mid\cF_{n})\\
&\hspace{-16em}=&\hspace{-8em}\bE_{x}(u(X_{n+1})\one_{\{\theta_2\leq n+1\}}\mid \cF_{n})
=\bE_{x}(\bE_{x}(u(X_{n+1})\one_{\{\theta_2\leq n+1\}}\mid
\mathcal{\widetilde{F}}_{n})\mid \mathcal{F}_{n})\nonumber\\
&\hspace{-16em}=&\hspace{-8em}\bE_{x}(\one_{\{\theta_2\leq n+1\}}\bE_{x}(u(X_{n+1})\mid
\mathcal{\widetilde{F}}_{n})\mid \mathcal{F}_{n})= \int_{\BbbE}u(y)f^2_{X_n}(y)\mu(dy)\bP_{x}(\theta_2 \leq n+1\mid \mathcal{F}_{n})\nonumber\\
&\hspace{-16em}\quad\stackrel{L.\ref{multidistform}}{=}&\hspace{-8em} \left(q_2\Pi^1_n+p_2\Pi_n^2+q_1\Pi^{12}_n
  \right)\int_{\BbbE}u(y)f^2_{X_n}(y)\mu(dy)  \nonumber
\end{eqnarray*}
We used the properties of conditional expectation and point 5 of Lemma~\ref{multidistform}. Similar transformations give us equations (\ref{wartOczPi1}), \eqref{wartOczPi12} and (\ref{wartOczPi1Pi2}) when the points 1 and 2, the point 4 and the point 1 of Lemma~\ref{multidistform}, respectively. From~\eqref{wartOczPi1}-\eqref{wartOczPi2} we get \eqref{exp1}. The proof of the lemma is complete.
\end{pf}

\section{Detection of new homogeneous segment}
\label{SolDetection1}
\subsection{Equivalent optimal stopping problem}
For $X_0 = x$ let us define: $Z_n = \bP_{x}(\theta_1 \leq n < \theta_2 \mid \mathcal{F}_n)$ for $n=0,1,2,\ldots$. We have
\begin{eqnarray}
\label{innyWzorNaZn}
    Z_n = \bP_{x}(\theta_1 \leq n < \theta_2 \mid \mathcal{F}_n) = \Pi_n^{1} - \Pi_n^2
\end{eqnarray}
$Y_n = \rm{ess}sup_{\{\tau \in \mathcal{T},\;\tau \geq n\}}\bP_{x}(\theta_1 \leq \tau < \theta_2 \mid \mathcal{F}_n)$ for $n=0,1,2,\ldots$ and
\begin{eqnarray}
\label{stopIntuicyjny}
    \tau_0 &=& \inf\{ n\geq 0: Z_n=Y_n \}
\end{eqnarray}

Notice that, if $Z_{\infty}=0$, then $Z_{\tau} = \bP_{x}(\theta_1 \leq \tau < \theta_2 \mid \mathcal{F}_{\tau})$ for
$\tau \in \mathcal{T}$. Since $\mathcal{F}_{n} \subseteq \mathcal{F}_{\tau}$ (when $n \leq \tau$) we have
\begin{eqnarray*}
    Y_n &=& \rm{ess}\sup_{\tau \geq n}E_{x}(Z_{\tau} \mid \mathcal{F}_n). \nonumber
\end{eqnarray*}
\begin{lemma}
    \label{lematCzasStopu}
The stopping time $\tau_0$ defined by the formula (\ref{stopIntuicyjny}) is the solution of the problem (\ref{problem}).
\end{lemma}

\begin{pf} From the theorems presented in \cite{boj79:disorder} it is enough to show that
$\displaystyle{\lim_{n \rightarrow \infty}Z_n=0}$. For all natural numbers $n,k$, where $n\geq k$ for each $x \in \mathbf{E}$ we have:
\begin{eqnarray*}
   Z_n &=& E_{x}(\one_{\{ \theta_1 \leq n < \theta_2 \}} \mid \mathcal{F}_n) \leq E_{x}(\sup_{j \geq n}\one_{\{ \theta_1 \leq j < \theta_2 \}} \mid \mathcal{F}_n)
\end{eqnarray*}
From Levy's theorem
$\limsup_{n\rightarrow \infty}Z_n \leq E_{x}(\sup_{j \geq k}\one_{\{ \theta_1 \leq j < \theta_2 \}} \mid \mathcal{F}_{\infty})$ where $\mathcal{F}_{\infty} = \sigma\left( \bigcup_{n=1}^{\infty}\mathcal{F}_n \right)$. It is true that: $\displaystyle{\lim_{k\rightarrow \infty}}\sup_{j \geq k}\one_{\{ \theta_1 \leq j < \theta_2 \}} = 0$ \emph{a.s.} and by the dominated convergence theorem we get
\[
  \lim_{k \rightarrow \infty}E_{x}(\sup_{j\geq k}\one_{\{ \theta_1 \leq j < \theta_2 \}} \mid \mathcal{F}_{\infty} ) = 0\;\; a.s.
\]
what ends the proof of the lemma.
\end{pf}

The reduction of the disorder problem to optimal stopping of Markov sequence is the consequence of the following lemma.
\begin{lemma}
    \label{RedukcjProblemu}
System $X^{x} = \left\{X^{x}_n\right\}$, where $X^{x}_n = (X_{n-1}, X_{n},\Pi_n^{1},\Pi_n^{2},\Pi_n^{12})$ forms a family of random Markov functions.
\end{lemma}
\begin{pf} Define a function:
\begin{equation}
\label{funkcjaFi}
    \varphi(x_1,x_2,\vecalpha\:;z) \\
    \quad= (x_2,z,\Pi^1(x_2,z,\vecalpha), \Pi^2(x_2,z,\vecalpha),\Pi^{12}(x_2,z,\vecalpha) )
\end{equation}

Observe that
\[
X^{x}_n = \varphi(X_{n-2}, X_{n-1},\vecPi_{n-1};X_n) = \varphi(X^{x}_{n-1};X_n)
\]
Hence $X^{x}_n$ can be interpreted as the function of the previous state $X^{x}_{n-1}$ and the random
variable $X_n$. Moreover, applying (\ref{exp1}), we get that the conditional distribution of $X_{n}$ given
$\sigma$-field $\mathcal{F}_{n-1}$ depends only on $X^{x}_{n-1}$. According to \cite{shi78:optimal} (pp. 102-103)
system $X^{x}$ is a family of random Markov functions.
\end{pf}
This fact implies that we can reduce the initial problem
(\ref{problem}) to the optimal stopping of the five-dimensional process
$(X_{n-1},X_{n},\Pi_n^{1},\Pi_{n}^{2},\Pi_n^{12})$ with the reward
\begin{equation}
\label{NowaWyplata}
        h(x_1,x_2,\vecalpha) = \alpha - \beta
\end{equation}
The reward function results from the equation (\ref{innyWzorNaZn}). Thanks to Lemma \ref{RedukcjProblemu} we construct the solution using standard tools of optimal stopping theory (cf \cite{shi78:optimal} ), as we do below.

Let us define two operators for any Borel function $v: \mathbf{E}^2\times [0,1]^{3}\longrightarrow [0,1]$ and the set $D=\{\omega:X_{n-1}=y, X_n=z,\Pi_n^{1} = \alpha, \Pi_{n}^{2} = \beta,\Pi_{n}^{12}=\gamma \}$:
\begin{eqnarray*}
T_{x}v(y,z,\vecalpha) &=& E_{x}(v(X_n, X_{n+1},\vecPi_{n+1})\mid D) \nonumber\\
\bQ_{x}v(y,z,\vecalpha) &=& \max\{v(y,z,\vecalpha), \bT_{x}v(y,z,\vecalpha) \} \nonumber
\end{eqnarray*}
From the well known theorems of optimal stopping theory (see \cite{shi78:optimal}), we infer that the solution of the problem \eqref{problem} is the Markov time $\tau_0$:
\begin{equation}\label{tau0star}
    \tau_0^\star = \inf\{n\geq 0: h(X_n, X_{n+1},\vecPi_{n+1})\geq h^{*}(X_n, X_{n+1},\vecPi_{n+1}) \},
\end{equation}
where:
\[
  h^{*}(y,z,\vecalpha) = \lim_{k\rightarrow \infty}\bQ_{x}^k h(y,z,\vecalpha).
\]
Of course
\[
    \bQ_{x}^k v(y,z,\vecalpha) = \max\{\bQ_{x}^{k-1} v, \bT_{x}\bQ_{x}^{k-1}v \} = \max\{ v, \bT_{x}\bQ_{x}^{k-1}v \}.
\]
To obtain a clearer formula for $\tau_0^\star$ and the solution of the problem \eqref{problem}, we formulate (cf \eqref{vecPi} and \eqref{vecalpha}):

\begin{theorem}
    \label{twierdz1}
\begin{description}
  \item[(a)]
 The solution \eqref{tau0star} of the optimal stopping problem for the stochastic system $X^x$ defined in Lemma~\ref{RedukcjProblemu} with payoff function (\ref{NowaWyplata}) is given by:
\begin{eqnarray}
\label{optymalnyStop}
    \tau_0^{*} = \inf\{n\geq 0:(X_n, X_{n+1},\vecPi_{n+1}) \in B^{*} \}.
\end{eqnarray}

Set $B^{*}$ is of the form:
\begin{eqnarray*}
B^{*} &=& \left\{ (y,z,\vecalpha) : (\alpha - \beta) \geq (1-\alpha-\gamma) \right. [ p_1\int_{\BbbE}R^{*}(y,u,\vecPi_1(y,u,\vecalpha))f^0_y(u)\mu(du) \nonumber \\
     &+& q_1\!\int_{\BbbE}S^{*}(y,u,\vecPi_1(y,u,\vecalpha))f^1_y(u)\mu(du) ] \nonumber \\
     &+& \left. (\alpha - \beta)p_2\!\int_{\BbbE}S^{*}(y,u,\vecPi_1(y,u,\vecalpha))f^1_y(u)\mu(du)\right\}, \nonumber
\end{eqnarray*}
where 
$ R^{*}(y,z,\vecalpha)= \lim_{k\rightarrow \infty}R^k(y,z,\vecalpha)$, $S^{*}(y,z,\vecalpha)= \lim_{k\rightarrow \infty}S^k(y,z,\vecalpha)$.
The functions $R^k$ and $S^k$ are defined recursively: $R^1(y,z,\vecalpha)=0$, $S^1(y,z,\vecalpha)=1$ and
\begin{eqnarray}
\label{R_k}
    &&R^{k+1}(y,z,\vecalpha) = (1-\one_{\calR_k}(y,z,\vecalpha))  \left(p_1\int_{\BbbE}R^{k}(y,u,\vecPi_1(y,u,\vecalpha))f^0_y(u)\mu(du)\right. \\
    &&\hspace{3em}\left. + q_1\!\int_{\BbbE}S^{k}(y,u,\vecPi_1(y,u,\vecalpha))f^1_y(u)\mu(du)\right),  \nonumber
\end{eqnarray}
\begin{eqnarray}\label{S_k}
S^{k+1}(y,z,\vecalpha) &=& \one_{\calR_k}(y,z,\vecalpha)+(1-\one_{\calR_k}(y,z,\vecalpha))\\
&&\quad\times p_2\!\int_{\BbbE}S^{k}(y,u,\vecPi_1(y,u,\vecalpha))f^1_y(u)\mu(du),\nonumber
\end{eqnarray}
where the set $\calR_k$ is:
\setlength\arraycolsep{2pt}
\begin{eqnarray}
\label{obszarR_k}
\calR_k &=& \left\{ (y,z,\vecalpha): h(y,z,\vecalpha) \geq \bT_{x}\bQ_{x}^{k-1}h(y,z,\vecalpha)  \right\} \\
     &=& \left\{ (y,z,\vecalpha):  (\alpha - \beta) \geq (1-\alpha-\gamma) \right.\nonumber\\
     &&\mbox{}\times \left[ p_1\int_{\BbbE}R^{k}(y,u,\vecPi_1(y,u,\vecalpha))f_y^0(u)\mu(du)\right.\nonumber\\
     &&\mbox{}+ \left.q_1\!\int_{\BbbE}S^{k}(y,u,\vecPi_1(y,u,\vecalpha))f^1_y(u)\mu(du) \right] \nonumber \\
     &&\mbox{}+ \left.(\alpha - \beta)p_2\! \int_{\BbbE}S^{k}(y,u,\vecPi^1(y,u,\vecalpha))f^1_y(u)\mu(du)\right\}.\nonumber
\end{eqnarray}
 \item[(b)] The optimal value for (\ref{problem}) is given by the formula
\setlength\arraycolsep{2pt}
\begin{align*}
     V(x)&=\max\{p_2\bar{\pi}\rho,V_0(x)\}\\
\intertext{where}
     V_0(x)&=\bar{\pi}\bar{\rho}\left[ p_1 \int_{\BbbE} R^{*}(x,u,\vecPi_1(x,u,\pi,\rho\pi,\rho\bar{\pi})) f^0_{x}(u)\mu(du)\right. \\
               &\left.+ q_1 \int_{\BbbE}S^{*}(x,u,\vecPi_1(x,u,\pi,\rho\pi,\rho(1-\pi)))f^1_{x}(u)\mu(du)\right]\\
               &+\bar{\pi}\rho p_2 \int_{\BbbE}S^{*}(x,u,\vecPi_1(x,u,\pi,\rho\pi,\rho(1-\pi)))f^1_{x}(u)\mu(du)
\end{align*}
and $\tau^\star=0\one_{\{p_2\bar{\pi}\rho\geq V_0(x)\}}+\tau_0^\star\one_{\{p_2\bar{\pi}\rho< V_0(x)\}}$.
\end{description}

\end{theorem}

\begin{pf} Part (a) results from Lemma~\ref{lematTechniczny} - the problem reduces to the optimal stopping of the Markov process $(X_{n-1}, X_{n},\Pi_n^{1},\Pi_{n}^{2},\Pi^{12}_n)$ with the payoff
function $h(y,z,\vecalpha) = \alpha - \beta$. Given (\ref{wartOczPi1Pi2}) with the function
$u$ equal to unity we get on $D=\{\omega:X_{n-1} = y, X_n=z,\Pi_n^1 = \alpha,\Pi_n^2=\beta,\Pi^{12}_n=\gamma \}$:
\begin{eqnarray*}
  \bT_{x}h(y,z,\vecalpha) &=&\bE_x\left(\Pi_{n+1}^{1}-\Pi_{n+1}^2 \mid \cF_n \right)\mid_{D} \nonumber\\
  &=& \left[((1-\Pi^{1}_n-\Pi^{12}_n)q_1+(\Pi^{1}_n-\Pi_n^2)p_2)\int_{\BbbE}f^1_{X_n}(u)\mu(du)\right] \mid_{D} \nonumber\\
  &=&(1-\alpha-\gamma)q_1 + (\alpha - \beta)p_2. \nonumber
\end{eqnarray*}
From the definition of $R^1$ and $S^1$ it is clear that
\[
    h(y,z,\vecalpha) = \alpha - \beta = (1- \alpha-\gamma)R^1(y,z,\vecalpha)+ (\alpha-\beta)S^1(y,z,\vecalpha)
\]
Also $\calR_1 = \{ (y,z,\vecalpha): h(y,z,\vecalpha) \geq \bT_{x}h(y,z,\vecalpha) \}$. 
From the definition of $\bQ_{x}$ and the facts above we obtain
\begin{eqnarray*}
    \bQ_{x}h(y,z,\vecalpha)&=& (1-\alpha-\gamma)R^2(y,z,\vecalpha) + (\alpha - \beta)S^2(y,z,\vecalpha), \nonumber
\end{eqnarray*}
where $ R^2(y,z,\vecalpha)=q_1(1-\one_{\calR_1}(y,z,\vecalpha))$ and
$S^2(y,z,\vecalpha)=p_2+(1-p_2)\one_{\calR_1}(y,z,\vecalpha))$.
Suppose the following induction hypothesis holds
\[
  \bQ_{x}^{k-1}h(y,z,\vecalpha) = (1-\alpha-\gamma)R^k(y,z,\vecalpha) + (\alpha - \beta)S^k(y,z,\vecalpha),
\]
where $R^k$ and $S^k$ are given by equations (\ref{R_k}), (\ref{S_k}), respectively. We will show 
\[
  \bQ_{x}^{k}h(y,z,\vecalpha) = (1-\alpha-\gamma)R^{k+1}(y,z,\vecalpha) + (\alpha - \beta)S^{k+1}(y,z,\vecalpha).
\]
From the induction assumption and the equations (\ref{wartOczPi1}), \eqref{wartOczPi12} and (\ref{wartOczPi1Pi2}) we obtain:
\setlength\arraycolsep{2pt}
\begin{eqnarray}
\label{operT}
\bT_{x}\bQ_{x}^{k-1}h(y,z,\vecalpha) &=& \bT_{x}(1-\alpha-\gamma)R^k(y,z,\vecalpha)\\
\nonumber&&\mbox{} + \bT_{x}(\alpha - \beta)S^k(y,z,\vecalpha)\\ 
&\hspace{-16em}=&\hspace{-8em} (1-\alpha-\gamma)p_1\int_{\BbbE}R^k(y,u,\vecPi_1(y,u,\vecalpha))f^0_y(u)\mu(du) \nonumber\\
   &\hspace{-16em}&\hspace{-8em}+ \left[(1-\alpha-\gamma)q_1+(\alpha-\beta)p_2\right]\int_{\BbbE}S^k(y,u,\vecPi_1(y,u,\vecalpha))f^1_y(u)\mu(du)\nonumber\\
\nonumber   &\hspace{-16em}=&\hspace{-8em} (1-\alpha-\gamma) \left[p_1\int_{\BbbE}R^k(y,u,\vecPi_1(y,u,\vecalpha))f^0_y(u)\mu(du)\right.\\
   &\hspace{-16em}&\hspace{-8em}\left. + q_1\int_{\BbbE}S^k(y,u,\vecPi_1(y,u,\vecalpha))f^1_y(u)\mu(du) \right] \nonumber\\     &\hspace{-16em}&\hspace{-8em}+(\alpha-\beta)p_2\int_{\BbbE}S^k(y,u,\vecPi_1(y,u,\vecalpha))f^1_y(u)\mu(du).\nonumber
\end{eqnarray}
Notice that
\begin{equation*}
(1-\alpha-\gamma)R^{k+1}(y,z,\vecalpha) + (\alpha - \beta)S^{k+1}(y,z,\vecalpha)
\end{equation*}
is equal $\alpha - \beta = h(y,z,\vecalpha) = \bQ_{x}^k h(y,z,\vecalpha)$ for $(y,z,\vecalpha) \in \mathcal{R}_k$ and, taking into account (\ref{operT}), it is equal $\bT_{x}\bQ_{x}^{k-1}h(y,z,\vecalpha) = \bQ_{x}^k h(y,z,\vecalpha)$ for $(y,z,\vecalpha) \notin \mathcal{R}_k$, where $\mathcal{R}_k$ is given by (\ref{obszarR_k}).
Finally we get
\[
  \bQ_{x}^{k}h(y,z,\vecalpha) = (1-\alpha-\gamma)R^{k+1}(y,z,\vecalpha) + (\alpha - \beta)S^{k+1}(y,z,\vecalpha).
\]
This proves (\ref{R_k}) and (\ref{S_k}). Using the monotone convergence theorem and the theorems of optimal stopping theory (see~\cite{shi78:optimal}) we conclude that the optimal stopping time $\tau_0^{*}$ is given by (\ref{optymalnyStop}).
\end{pf}

\begin{pf} Part (b). First, notice that $\Pi_1^{1}$, $\Pi^2_1$ and $\Pi^{12}_1$ are given by \eqref{eqpi1}-\eqref{eqpi12} 
and the boundary condition formulated in Lemma~\ref{reqform}.
Under the assumption $\tau^{*} < \infty$ a.s. we get:
\setlength\arraycolsep{2pt}
\begin{eqnarray*}
&&\bP_{x}(\tau^{*} < \infty, \theta_1 \leq \tau^{*}\!<\theta_2) = \sup_{\tau} \bE Z_{\tau} \nonumber \\
&&\quad= \bE\max\{ h(x, X_{1},\vecPi_{1}), \bT_{x}h^{*}(x, X_{1},\vecPi_{1}) \} = \bE\lim_{k\rightarrow \infty}\!\bQ_{x}^k h(x, X_{1},\vecPi_{1})\nonumber\\
&&\quad= \bE\left[ (1-\Pi_1^{1}-\Pi^{12}_1)R^{*}(x, X_{1},\vecPi_{1}) +(\Pi_1^{1}-\Pi_1^2)S^{*}(x, X_{1},\vecPi_{1}) \right] \nonumber\\
&&\quad= \bar{\pi}\bar{\rho}p_1\!\int_{\BbbE}R^{*}(x,u,\vecPi_1(x,u,\pi,\rho\pi,\rho\bar{\pi}))f^0_x(u)\mu(du)\\
&&+(\bar{\pi}\bar{\rho} q_1+\pi\bar{\rho}p_2)\!\int_{\BbbE}S^{*}(x,u,\vecPi_1(x,u,\pi,\rho\pi,\rho\bar{\pi}))f^1_{x}(u)\mu(du).\nonumber
\end{eqnarray*}
\normalsize
We used Lemma \ref{lematTechniczny} here and simple calculations for $\Pi_{1}^{1}$, $\Pi_{1}^{2}$ and $\Pi^{12}_1$.
This ends the proof.
\end{pf}

\subsection{Remarks}
    It is notable that the solution of formulated problem depends only on two-dimensional vector of
posterior processes because $\Pi_n^{12}=\rho(1-\Pi^1_n)$. The obtained formulae  are very general and for this reason -- quite complicated. We simplify the model by assuming that $P(\theta_1>0)=1$ and $P(\theta_2>\theta_1)=1$. However, it seems that some further simplifications can be made in special cases. Further research should be carried out in this direction. From a practical point of view, computer algorithms are necessary to construct $B^{*}$ -- the set in which it is optimally to stop our observable sequence.

\section{Immediate detection of the first and the second disorder}
\label{SolDisorders1}
\subsection{Equivalent double optimal stopping problem}\label{dosp}
Let us consider the problem $\mbox{D}_{00}$ formulated in \eqref{equ3}. A {\it compound stopping variable} is a pair $(\tau,\sigma)$ of stopping times such that $0\leq\tau\leq\sigma$ a.e.. The aim is to find a compund stopping variable $(\tau^\star,\sigma^\star)$ such that
\small
\begin{equation}
\label{equ3A}
\bP_x((\theta_1,\theta_2)=(\tau^*,\sigma^*)) =
\sup_{\stackrel{(\tau,\sigma)\in\cT}{0\leq\tau\leq\sigma<\infty}} \bP_x((\theta_1,\theta_2)=(\tau,\sigma)).
\end{equation}
\normalsize
Denote $\cT_m=\{(\tau,\sigma)\in\cT: \tau\geq m\}$, 
$\cT_{mn}=\{(\tau,\sigma)\in\cT: \tau=m, \sigma\geq n\}$ and $\cS_m=\{\tau\in\cS:\tau\geq m\}$. Let us denote
$\cF_{mn}=\cF_n$, $m,n\in\BbbN$, $m\leq n$. We define two-parameter stochastic sequence 
$\xi(x)=\{\xi_{mn},\ m,n\in\BbbN,\  m<n,\ x\in\BbbE\}$, where
\[
\xi_{mn}=\bP_{x}(\theta_1=m, \theta_2=n|\cF_{mn}).
\]
We can consider for every $x\in\BbbE$, $m,n\in\BbbN$, $m<n$, the optimal stopping problem of $\xi(x)$ on 
$\cT^\plus_{mn}=\{(\tau,\sigma)\in\cT_{mn}:\tau<\sigma\}$.  
A compound stopping variable $(\tau^*,\sigma^*)$ is said to be 
optimal in $\cT^\plus_m$ (or $\cT^\plus_{mn}$) if 
\begin{equation}\label{doubleDISmarkovian}
\bE_{x}\xi_{\tau^*\sigma^*} =\sup_{(\tau,\sigma)\in\cT_m}\bE_{x}\xi_{\tau\sigma} 
\end{equation}
(or $\bE_{x}\xi_{\tau^*\sigma^*} =\sup_{(\tau,\sigma)\in\cT^\plus_{mn}}\bE_{x}\xi_{\tau\sigma}$). 
Let us define
\begin{equation}
\label{pr1}
\eta_{mn} =
\esssup_{(\tau,\sigma)\in\cT^\plus_{mn}}\bE_x(\xi_{\tau\sigma}|\cF_{mn}).
\end{equation}
If we put $\xi_{m\infty}=0$, then 
\[
\eta_{mn}=\esssup_{(\tau,\sigma)\in\cT^\plus_{mn}}\bP_{x}(\theta_1=\tau,\theta_2=\sigma|\cF_{mn}).
\]
From the theory of optimal stopping for double indexed processes (cf. \cite{hag67:more},\cite{nik81}) the sequence $\eta_{mn}$
satisfies 
\[
\eta_{mn}=\max\{\xi_{mn},\bE(\eta_{mn+1}|\cF_{mn})\}.
\]
Moreover, if $\sigma^*_m=\inf\{n>m: \eta_{mn}=\xi_{mn}\}$, then
$(m,\sigma^*_n)$ is optimal in $\cT^\plus_{mn}$ and $\eta_{mn}=\bE_x(\xi_{m\sigma^*_n}|\cF_{mn})$ a.e.. The case when there are no segment with the distribution $f^1_x(y)$ appears with probability $\rho$. It will be taken into account. Define
\[
\hat{\eta}_{mn}=\max\{\xi_{mn},\bE(\eta_{m\;n+1}|\cF_{mn})\}, \text{ for $n\geq m$.}
\]
if $\hat{\sigma}^*_m=\inf\{n\geq m: \hat{\eta}_{mn}=\xi_{mn}\}$, then $(m,\hat{\sigma}^*_m)$ is optimal in $\cT_{mn}$ and $\hat{\eta}_{mm}=\bE_x(\xi_{m\sigma^*_m}|\cF_{mm})$ a.e.. For further consideration denote
\begin{equation}
\label{pr2}
\eta_{m} = \bE_{x}(\eta_{mm+1}|\cF_m).
\end{equation} 
\begin{lemma}
\label{lem1}
The stopping time $\sigma^*_m$ is optimal for every stopping problem
(\ref{pr1}).
\end{lemma}
\begin{pf}
It suffices to prove $\lim_{n\rightarrow\infty}\xi_{mn}=0$ (cf.
\cite{boj79:disorder}). We have for $m,n,k\in\BbbN$, $n\geq k>m$ and every
$x\in\BbbE$ 
\[
\bE_x(\one_{\{\theta_1=m,\theta_2=n\}}|\cF_{mn})=\xi_{mn}(x)\leq
\bE_x(\sup_{j\geq k}\one_{\{\theta_1=m,\theta_2=j \}}|\cF_{m}),
\]
where $\one_A$ is the characteristic function of the set $A$.
By Levy's theorem
\[
\limsup_{n\rightarrow\infty}\xi_{mn}(x)\leq\bE_x(\sup_{j\geq
k}\one_{\{\theta_1=m,\theta_2=j \}}|\cF_{n\infty}),
\]
where $\cF_\infty=\cF_{n\infty}=\sigma(\bigcup_{n=1}^\infty
\cF_{n})$.  We have
$\displaystyle{\lim_{k\rightarrow\infty}\sup_{j\geq k}}\mbox{}\one_{\{\theta_1=m,\theta_2=j\}}=0 \mbox{ a.e.}$ and by dominated convergence theorem 
\[
\lim_{k\rightarrow\infty}\bE_x(\sup_{j\geq k}
\mbox{}\one_{\{\theta_1=m,\theta_2=j\}} | \cF_\infty)=0.
\]
\end{pf}

What is left is to consider the optimal stopping problem for ${(\eta_{mn})_{m=0,}^{\infty,}}_{n=m}^{\infty}$ on $(\cT_{mn})_{m=0,n=m}^{\infty,\infty}$. 
Let us define 
\begin{equation}
\label{FirstStopPayment}
V_m=\esssup_{\tau\in\cS_m}\bE_x(\eta_\tau|\cF_m).
\end{equation}
Then $V_m=\max\{\eta_m,\bE_x(V_{m+1}|\cF_m)\} \mbox{ a.e.}$ and
we define $\tau^*_n=\inf\{k\geq n: V_k=\eta_k\}$. 
\begin{lemma}
\label{lem2}
The strategy $\tau^*_0$ is the optimal first stop.
\end{lemma}
\begin{pf}
To show that $\tau^{*}_0$ is the optimal first stop strategy we prove that
$\bP_x(\tau^*_0<\infty)=1$. To this end, we argue in the usual
manner i.e. we show $\lim_{m\rightarrow\infty} \eta_m=0$.

We have 
\begin{eqnarray*}
\eta_m&=&\bE_x(\xi_{m\sigma_m^*}|\cF_m)
=\bE_x(\bE_x(\one_{\{\theta_1=m,\theta_2 =\sigma_m^*\}}|\cF_{m\sigma_m^*})|\cF_m) \\
 &=&\bE_x(\one_{\{\theta_1=m,\theta_2=\sigma_m^*\}}|\cF_m)
\leq\bE_x(\sup_{j\geq k}\one_{\{\theta_1=j,\theta_2=\sigma^*_j\}}|\cF_m).
 \end{eqnarray*}
Similarly as in proof of Lemma \ref{lem1} we have got
 \[
\limsup_{m\rightarrow\infty}\eta_m(x)\leq\bE_x(\sup_{j\geq k}
\mbox{}\one_{\{\theta_1=j,\theta_2=\sigma^*_j\}}|\cF_\infty). 
\]
Since 
$\lim_{k\rightarrow\infty}\sup_{j\geq k}\one_{\{\theta_1=k,\theta_2=\sigma^*_j\}} \leq 
\limsup_{k\rightarrow\infty}\one_{\{\theta_1=k\}}=0$, 
it follows that
\[
\lim_{m\rightarrow\infty}\eta_m(x)\leq\lim_{k\rightarrow\infty}
\bE_x(\sup_{j\geq k} \one_{\{\theta_1=j,\theta_2=\sigma_j^*\}}|
\cF_\infty)=0.  
\]
\end{pf}

Lemmas \ref{lem1} and \ref{lem2} describe the method of solving
the ``disorder problem'' formulated in Section \ref{sformProblem} (see \eqref{equ3A}).

\subsection{Solution of the equivalent double stopping problem}
\label{sds4}

For the sake of simplicity we shall confine ourselves to the case $d_1=d_2=0$. It will be easily seen how to generalize the solution of the problem to solve $\mbox{D}_{d_1d_2}$ for $d_1>0$ or $d_2>0$. First of all we construct multidimensional  Markov chains such that $\xi_{mn}$ and $\eta_m$ will be the functions of their states. By consideration of Section~\ref{PosterioriProc} concerning \emph{a posteriori} processes we get $\xi_{00}=\pi\rho$ and for $m<n$
\begin{eqnarray*}
\xi_{m\,n}^x &=& \bP_x(\theta_1=m,\theta_2=n|\cF_{m\,n})\\
 &=&\bar{\pi}\bar{\rho}\frac{p_1^{m-1}q_1p_2^{n-m-1}q_2 
\prod_{s=1}^{j-1}f^0_{X_{s-1}}(X_s)\prod_{t=j}^{n-1}f^1_{X_{t-1}}(X_{t})
f^2_{X_{n-1}}(X_n)}{S_n(x_0,X_1,\ldots,X_n)}\\
 &=& \frac{q_2}{p_2}\Pi_{m\,n}(x)\frac{f^2_{X_{n-1}}(X_n)}{f^1_{X_{n-1}}(X_n)}
\end{eqnarray*}
and for $n=m$, by Lemma \ref{multidistform}, 
\begin{equation}
\label{ximm}
\xi_{m\,m}^x = \bP_x(\theta_1=m,\theta_2=m|\cF_{m\,m})=\rho\frac{q_1}{p_1}\frac{f^2_{X_{m-1}}(X_m)}{f^0_{X_{m-1}}(X_m)}(1-\Pi_m^1).
\end{equation}
We can observe that $(X_n,X_{n+1},\vecPi_{n+1},\Pi_{m\,n+1})$ for $n=m+1, m+2,\ldots$ is a function of
$(X_{n-1},X_{n},\vecPi_{n},\Pi_{m\,{n}})$ and $X_{n+1}$. Besides, the conditional distribution of $X_{n+1}$
given $\cF_n$ (cf. (\ref{exp1})) depends on $X_n$, $\Pi^1_n(x)$
and $\Pi^2_n(x)$ only. These facts imply that  $\{(X_n,X_{n+1},\vecPi_{n+1},\Pi_{m\,n+1})\}_{n=m+1}^\infty$
form a homogeneous Markov process (see Chapter~2.15 of \cite{shi78:optimal}). This allows us to reduce the problem
(\ref{pr1}) for each $m$ to the optimal stopping problem of the Markov process
$Z_m(x)=\{(X_{n-1},X_n,\vecPi_n,\Pi_{m\,n}),\ m,n\in\BbbN,\ 
 m<n,\ x\in\BbbE\}$ with the reward function $h(t,u,\vecalpha,\delta)= \frac{q_2}{p_2}\delta\frac{f^2_t(u)}{f^1_t(u)}$. 
 
\begin{lemma}
 \label{lem4}
 A solution of the optimal stopping problem (\ref{pr1}) for $m=1,2,\ldots$ has a form
\begin{equation}
\label{optsstop}
\sigma^*_m=\inf\{n>m:
\frac{f^2_{X_{n-1}}(X_n)}{f^1_{X_{n-1}}(X_n)}\geq R^*(X_n)\}
\end{equation}
where $R^*(t)=p_2\int_{\BbbE}r^*(t,s)f^1_t(s)\mu(ds)$. The function
$r^*=\lim_{n\rightarrow\infty} r_n$, where
$r_0(t,u)=\frac{f^2_{t}(u)}{f^1_{t}(u)}$, 
\begin{equation}
\label{rmale}
r_{n+1}(t,u)=\max\{\frac{f^2_t(u)}{f^1_t(u)},p_2\int_{\BbbE}r_n(u,s)f^1_u(s)\mu(ds)\}.
\end{equation}
So $r^*(t,u)$ satisfies the equation
\begin{equation}
r^*(t,u)=\max\{\frac{f^2_t(u)}{f^1_t(u)},p_2\int_{\BbbE}
r^*(u,s)f^1_u(s)\mu(ds)\}.
\end{equation}
The value of the problem 
\begin{equation}
\label{valpr1}
\eta_m = \bE_x(\eta_{m\,{m+1}}|\cF_m)
      =\frac{q_1}{p_1}\frac{f^1_{X_{m-1}}(X_m)}{f^0_{X_{m-1}}(X_m)}(1-\Pi^1_m)R_\rho^\star(X_{m-1},X_m),
\end{equation}
where
\begin{equation}
\label{Rrho}       
R_\rho^\star(t,u)=\max\{\rho\frac{f^2_{t}(u)}{f^1_{t}(u)},\frac{q_2}{p_2}(1-\rho)R^\star(u)\}.
\end{equation}
\end{lemma}

\begin{pf}
For any Borel function $u:\BbbE\times\BbbE\times [0,1]^4\rightarrow
[0,1]$ and $D=\{\omega:X_{n-1}=t,X_n=u,\Pi^1_n(x)=\alpha,\Pi^2_n(x)=\beta,\Pi^{12}_n=\gamma,\Pi_{m\,n}(x)=\delta\}$ let us define two operators
\begin{eqnarray*}
\bT_xu(t,u,\vecalpha,\delta) &=&
\bE_x(u(X_n,X_{n+1},\vecPi_{n+1}(x),\Pi_{m\,n+1}(x))|D) 
\end{eqnarray*}
and 
\[
\bQ_xu(t,u,\vecalpha,\delta)=\max\{u(t,u,\vecalpha,\delta),\bT_xu(t,u,\vecalpha,\delta)\}.
\]
On the bases of the well-known theorem from the theory of optimal stopping (see \cite{shi78:optimal}, \cite{nik81})
we conclude that the solution of (\ref{pr1}) is a Markov time
\[
\sigma^*_m=\inf\{n>m: h(X_{n-1},X_n,\vecPi_n,\Pi_{m\,n})=h^*(X_{n-1},X_n,\vecPi_n(x),\Pi_{m\,n})\},
\]
where $h^*=\lim_{k\rightarrow\infty}\bQ^k_xh(t,u,\vecalpha,\delta)$.
By \eqref{eqpi} and (\ref{exp1}) on $D=\{\omega:X_{n-1}=t,X_{n}=u,\Pi^1_n=\alpha,\Pi^2_n=\beta,\Pi^{12}_n=\gamma,\Pi_{m\,n}=\delta\}$ we have
\begin{eqnarray*}
\bT_xh(t,u,\vecalpha,\delta)&=&\bE_x(\frac{q_2}{p_2}\Pi_{m\,n+1}
\frac{f^2_{X_n}(X_{n+1})}{f^1_{X_n}(X_{n+1})}|D)\\
		 &=&\frac{q_2}{p_2}\delta p_2\bE(\frac{f^1_{u}(X_{n+1})}{H(u,X_{n+1},\vecalpha)}
		 \frac{f^2_u(X_{n+1})}{f^1_u(X_{n+1})}|\cF_n)|{}_D\\
&\stackrel{(\ref{exp1})}{=} &q_2\delta\int_{\BbbE}\frac{f^2_u(s)}{H(u,s,\vecalpha)} 
H(u,s,\vecalpha)\mu(ds) = q_2\delta
\end{eqnarray*}
and
\begin{equation}
\label{max1}
\bQ_xh(t,u,\vecalpha,\delta) =
\frac{q_2}{p_2}\delta\max\{\frac{f^2_t(u)}{f^1_t(u)},p_2\}.
\end{equation}
Let us define $r_0(t,u)=1$ and
\[
r_{n+1}(t,u) = \max\{\frac{f^2_t(u)} {f^1_t(u)} ,
p_2\int_{\BbbE}r_n(u,s)f^1_u(s)\mu(ds)\}.
\]
We show that 
\begin{equation}
\label{maxn}
\bQ_x^{\ell} h(t,u,\vecalpha,\delta)=\frac{q_2}{p_2}\delta r_{\ell}(t,u)
\end{equation}
for $\ell = 1,2,\ldots$. 
We have by (\ref{max1}) that $\bQ_xh=\frac{q_2}{p_2}\gamma r_1$. Let us assume (\ref{maxn}) for $\ell\leq k$. By (\ref{exp1}) on $D=\{\omega:X_{n-1}= t,X_{n}=u,\Pi^1_n=\alpha,\Pi^2_n=\beta,\Pi^{12}_n=\gamma,\Pi_{mn}=\delta\}$ we have got
\begin{eqnarray*}
\bT_x\bQ_x^kh(t,u,\vecalpha,\delta)&=&\bE_x(\frac{q_2}{p_2}\Pi_{m\,k+1}
r_k(X_n,X_{n+1})|D) \\
 &=& \frac{q_2}{p_2}\delta p_2\int_{\BbbE}r_k(u,s)f^1_u(s)\mu(ds).
\end{eqnarray*}
It is easy to show (see \cite{shi78:optimal}) that
\[
\bQ^{k+1}_xh=\max\{h,\bT_x\bQ^k_xh\}, \mbox{ for $k=1,2,\ldots$}.
\]
Hence we have got
$\bQ^{k+1}_xh=\frac{q_2}{p_2}\delta r_{k+1}$ and (\ref{maxn}) is
proved for $\ell = 1,2,\ldots$. This gives 
\begin{equation}
\label{hstar}
h^*(t,u,\vecalpha,\delta)=\frac{q_2}{p_2}\delta\lim_{k\rightarrow\infty}
r_{k}(t,u)=\frac{q_2}{p_2}\delta r^*(t,u)
\end{equation}
and
\[
\eta_{m\,n} = \esssup_{(\tau,\sigma)\in\cT_{m\,n}}
\bE_x(\xi_{\tau,\sigma} |\cF_{m\,n})=h^*(X_{n-1},X_{n},\vecPi_n,\Pi_{m\,n}).
\]
We have  by \eqref{hstar} and \eqref{eqpi}
\[
\bT_xh^*(t,u,\vecalpha,\delta)=\frac{q_2}{p_2}\delta
p_2\int_{\BbbE}r^*(u,s)f^1_u(s)\mu(ds)=\frac{q_2}{p_2}\delta R^*(u)
\]
and $\sigma^*_m$ has form (\ref{optsstop}).
By (\ref{pr2}), \eqref{ximm} and (\ref{exp1}) we obtain
\small
\begin{eqnarray}
\label{payoffForFirst}
\eta_m&=&\max\{\xi^x_{mm},\bE(\eta_{m\,m+1}|\cF_m)\}=f(X_{m-1},X_m,\vecPi_m,\Pi_{mm})\\
\nonumber&=&\max\{\rho\frac{q_1}{p_1}\frac{f^2_{X_{m-1}}(X_m)}{f^0_{X_{m-1}}(X_m)}(1-\Pi^1_m),\frac{q_2}{p_2}(1-\Pi_{mm})R^\star(X_m)\}\\
\nonumber&\stackrel{L.\ref{reqform}}{=}&\frac{q_1}{p_1}\frac{f^1_{X_{m-1}}(X_m)}{f^0_{X_{m-1}}(X_m)}(1-\Pi^1_m)R_\rho^\star(X_{m-1},X_m).
\end{eqnarray}
\normalsize
\end{pf}
\begin{remark}
Based on the results of Lemma~\ref{lem4} and properties of the \emph{a posteriori} process $\Pi_{nm}$ we have that the expected value of success for the second stop when the observer stops immediately at $n=0$ is $\pi\rho$ and when at least one observation has been made $\bE(\eta_1|\cF_0)=\frac{q_1}{p_1}\bE((1-\Pi_1^1)\frac{f_x^1(X_1)}{f_x^0(X_1)}R^\star_\rho(x,X_1)|\cF_0)=
\frac{q_1}{p_1}(1-\pi)p_1\int_\BbbE f_x^1(u)R^\star_\rho(x,u)\mu(du)$. As a consequence we have optimal second moment 
\[
\hat\sigma_0^\star=\left\{\begin{array}{ll}
0&\mbox{if $\pi\rho\geq q_1 (1-\pi)\int_\BbbE f_x^1(u)R^\star_\rho(x,u)\mu(du)$,}\\
\sigma_0^\star&\mbox{ otherwise.}
\end{array}
\right.
\]
\end{remark}
By Lemmas \ref{lem4} and \ref{reqform} (the formula \eqref{eqpi}) the optimal stopping problem
\eqref{FirstStopPayment} has been transformed to the optimal stopping
problem for the homogeneous Markov process
$$W=\{(X_{m-1},X_m,\vecPi_m),\ m\in\BbbN,\ x\in\BbbE\}$$
with the reward function
\begin{equation}
\label{FstepPayoff}
f(t,u,\vecalpha)=\frac{q_1}{p_1}\frac{f^1_{t}(u)}{f^0_{t}(u)}(1-\alpha) R^\star_\rho(t,u).
\end{equation}

\begin{theorem}
 \label{lem5}
 A solution of the optimal stopping problem \eqref{FirstStopPayment} for $n=1,2,\ldots$ has a form
\begin{equation}
\label{optfstop}
\tau^*_n=\inf\{k\geq n: (X_{k-1},X_k,\vecPi_k,)\in B^*\}
\end{equation}
where $B^*=\{(t,u,\vecalpha): \frac{f^2_t(u)}{f^1_t(u)}R^\star_\rho(t,u)\geq
p_1\int_{\BbbE} v^*(u,s)f^0_u(s)\mu(ds)\}$. The function
$v^*(t,u)=\lim_{n\rightarrow\infty} v_n(t,u)$, where
$v_0(t,u)=R^\star_\rho(t,u)$, 
\begin{equation}
\label{vmale}
v_{n+1}(t,u)=\max\{\frac{f^2_t(u)}{f^1_t(u)}R^\star_\rho(t,u),p_1\int_{\BbbE}
v_n(u,s)f^1_u(s)\mu(ds)\}.
\end{equation}
So $v^*(t,u)$ satisfies the equation
\begin{equation}
v^*(t,u)=\max\{\frac{f^2_t(u)}{f^1_t(u)}R^\star_\rho(t,u),p_1\int_{\BbbE}v^*(u,s)f^1_u(s)\mu(ds)\}.
\end{equation}
The value of the problem $V_n=v^*(X_{n-1},X_n)$.
\end{theorem}

\begin{pf}
For any Borel function $u:\BbbE\times\BbbE\times [0,1]^3\rightarrow [0,1]$ and $D=\{\omega:X_{n-1}=t,X_n=u,\Pi^1_n(x)=\alpha,\Pi^2_n(x)=\beta,\Pi^{12}_n=\gamma\}$ let us define two operators
\begin{eqnarray*}
\bT_xu(t,u,\vecalpha) &=& \bE_x(u(X_n,X_{n+1},\vecPi_{n+1})|D ) 
\end{eqnarray*}
and $\bQ_xu(t,u,\vecalpha)=\max\{u(t,u,\vecalpha),\bT_xu(t,u,\vecalpha)\}$.
Similarly as in the proof of Lemma \ref{lem4} 
we conclude that the solution of \eqref{FirstStopPayment} is a Markov time 
\[
\tau^*_m=\inf\{n>m: f(X_{n-1},X_n,\vecPi_n)=f^*(X_{n-1},X_n,\vecPi_n)\},
\]
where $f^*=\lim_{k\rightarrow\infty}\bQ^k_xf(t,u,\vecalpha)$.
By (\ref{exp1}) and \eqref{FstepPayoff} on $D=\{\omega:X_{n-1}=t,X_{n}=u,\Pi^1_n=\alpha,\Pi^2_n=\beta,\Pi^{12}_n=\gamma\}$ we have
\begin{eqnarray*}
\bT_xf(t,u,\vecalpha)&=&\bE_x(\frac{q_1}{p_1}(1-\Pi^1_{n+1})
\frac{f^1_{X_n}(X_{n+1})}{f^0_{X_n}(X_{n+1})}R_\rho^\star(X_n,X_{n+1})|D)\\
		 &=&\frac{q_1}{p_1}(1-\alpha)p_1\bE(\frac{f^0_{u}(X_{n+1})}{H(u,X_{n+1},\alpha,\beta)} 
		 \frac{f^1_u(X_{n+1})}{f^0_u(X_{n+1})}R_\rho^\star(X_n,X_{n+1})|\cF_n)|{}_D\\
&\stackrel{(\ref{exp1})}{=} &
\frac{q_1 }{p_1 }(1-\alpha) p_1\int_{\BbbE}\frac{f^1_u(s)}{H(u,s,\alpha,\beta)}  
H(u,s,\alpha,\beta)R_\rho^*(u,s)\mu(ds) \\
&=&\frac{q_1 }{p_1}(1-\alpha) p_1\int_{\BbbE}R_\rho^*(u,s)f^1_{u}(s) \mu(ds)
\end{eqnarray*}
and
\begin{eqnarray}
\label{max2}
\bQ_xf(t,u,\vecalpha) &=& \frac{q_1}{p_1 } (1-\alpha)\max \{\frac{ f^1_t(u)} {f^0_t(u)}R^\star_\rho(t,u),
p_1\int_{\BbbE}R_\rho^\star(u,s)f^1_{u}(s)\mu(ds)\}\\
\nonumber &=&\frac{q_1}{p_1} \alpha v_1(t,u).
\end{eqnarray}
Let us define $v_1(t,u)=\max \{\frac{ f^1_t(u)}{f^0_t(u)}R^\star_\rho(t,u),p_1\int_{\BbbE}R_\rho^*(u,s)f^1_{u}(s)\mu(ds)$ and
\[
v_{n+1}(t,u) = \max\{\frac{f^1_t(u)} {f^0_t(u)}R^\star_\rho(t,u),p_1\int_{\BbbE}v_n(u,s)f^0_u(s)\mu(ds)\}.
\]
We show that 
\begin{equation}
\label{maxn2}
\bQ_x^{\ell} f(t,u,\vecalpha)=\frac{q_1}{p_1 }(1-\alpha) v_{\ell}(t,u)
\end{equation}
for $\ell = 1,2,\ldots$. 
We have by (\ref{max2}) that $\bQ_xf(t,u,\vecalpha)=\frac{q_1}{p_1}(1-\alpha) v_1(t,u)$ and let us 
assume (\ref{maxn2}) for $\ell\leq k$. By (\ref{exp1}) on $D=\{\omega:X_{n-1} = t ,
X_{n}=u,\Pi^1_n=\alpha,\Pi^2_n=\beta,\Pi^{12}_n=\gamma\}$ we have got
\begin{eqnarray*}
\bT_x\bQ_x^kf(t,u,\vecalpha)&=&\bE_x(\frac{q_1}{p_1}(1-\Pi^1_{k+1})v_k(X_n,X_{n+1})|D) \\
 &=& \frac{q_1}{p_1}(1-\alpha) p_1\int_{\BbbE}v_k(u,s)f^0_u(s)\mu(ds).
\end{eqnarray*}
Hence we have got
$\bQ^{k+1}_xf=\frac{q_1}{p_1}(1-\alpha) v_{k+1}$ and (\ref{maxn2}) is
proved for $\ell = 1,2,\ldots$. This gives 
\[
f^*(t,u,\vecalpha)=\frac{q_1}{p_1}(1-\alpha)\lim_{k\rightarrow\infty}
v_{k}(t,u)=\frac{q_1}{p_1}\alpha v^*(t,u)
\]
and
\[
V_{m} = \frac{q_1}{p_1}(1-\Pi^1_m) v^*(X_{m-1},X_m).
\]
We have 
\[
\bT_xf^*(t,u,\vecalpha)=\frac{q_1}{p_1}(1-\alpha)
p_1\int_{\BbbE}v^*(u,s)f^0_u(s)\mu(ds).
\]
Define $B^*= \{(t,u,\vecalpha): \frac{f^1_t(u)}{f^0_t(u)}R_\rho^\star(t,u)\geq
p_1\int_{\BbbE} v^*(u,s)f^0_u(s)\mu(ds)\}$ then $\tau^*_n$ for $n\geq 1$ has a form
(\ref{optfstop}). The value of the problem \eqref{doubleDISmarkovian}, \eqref{FirstStopPayment} and \eqref{equ3} is equal
\[
v_0(x)=\max\{\pi,\bE_x(V_1|\cF_0)\}=\max\{\pi,\frac{q_1}{p_1}(1-\pi)p_1\int_{\BbbE}v^*(u,s)f^0_u(s)\mu(ds)\}
\]
and 
\[
\hat\tau_0^*=\left\{\begin{array}{ll}
0&\mbox{ if $\pi\geq q_1(1-\pi)\int_{\BbbE}v^*(u,s)f^0_u(s)\mu(ds)$,}\\
\tau^*_0&\mbox{otherwise.}
\end{array}\right.
\]
\end{pf}		

Based on Lemmas \ref{lem4} and \ref{lem5} the solution of the problem
$\mbox{D}_{00}$ can be formulated as follows.

\begin{theorem}
\label{th1}
A compound stopping time $(\tau^*,\sigma^*_{\tau^*})$, where
$\sigma^*_m$ is given by (\ref{optsstop}) and $\tau^*=\hat\tau^*_0$ is given by
(\ref{optfstop}), is the solution of the problem $\mbox{D}_{00}$. The value of the problem
\[
\bP_x(\tau^*<\sigma^*<\infty,\theta_1=\tau^*,\theta_2=\sigma^*_{\tau^*})=
\max\{\pi,q_1(1-\pi)\int_{\BbbE}v^*(u,s)f^0_u(s)\mu(ds)\}.
\]
\end{theorem}

\begin{remark}
The problem can be extended to optimal detection of more than two successive disorders. The distribution of $\theta_1$, $\theta_2$ may be more general. The general {\it a priori} distributions of disorder moments leads to more complicated formulae, since the corresponding Markov chains are not homogeneous. 
\end{remark}

\section{Final remarks}
    It is notable that the final optimal solutions turns out to have an unexpectedly simple form. It seems that some further simplifications can be made in special cases. From a practical point of view, computer algorithms are necessary to construct $B^{*}$ -- the set in which we stop our observable sequence. Since we always refer to the transitions densities it is still open problem of switching between the independent Markov sequences. 

\appendix{Useful relations}
\subsection{\label{distDISsample}Distributions of disordered samples}
Let us introduce the $n$-dimensional distribution for various configuration of disorders.
\begin{eqnarray}
\label{theta1lesstheta2lessn} f^{\theta_1\leq\theta_2\leq n}_x(\vec{x}_{1,n}) &=&\bar{\pi}\rho\sum_{j=1}^{n}\{p_1^{j-1}q_1
\prod_{s=1}^{j-1}f^0_{x_{s-1}}(x_s)\prod_{t=j}^{n}f^2_{x_{t-1}}(x_{t})\}\\
\nonumber &\hspace{-12em}+&\hspace{-7em}\bar{\pi}\bar{\rho} \sum_{j=1}^{n-1}\sum_{k=j+1}^{n}\{p_1^{j-1}q_1p_2^{k-j-1}q_2
\prod_{s=1}^{j-1}f^0_{x_{s-1}}(x_s)\prod_{t=j}^{k-1}f^1_{x_{t-1}}(x_{t}) \prod_{u=k}^{n}f^2_{x_{u-1}}(x_u)\}\\
\nonumber&\hspace{-12em}+&\hspace{-7em}\pi\rho \prod_{s=1}^n f^2_{x_{s-1}}(x_s)\\
\label{theta1less_n_lesstheta2} f^{\theta_1\leq n<\theta_2}_x(\vec{x}_{1,n}) &=&\bar{\pi}\bar{\rho}\sum_{j=1}^n\{p_1^{j-1}q_1 p_2^{n-j}\prod_{s=1}^{j-1}f^0_{x_{s-1}}(x_s)\prod_{t=j}^{n}f^1_{x_{t-1}}(x_t)\}\\
\nonumber &+&\pi\bar{\rho}\sum_{j=1}^n\{p_2^{j-1}q_2 \prod_{s=1}^{j-1}f^1_{x_{s-1}}(x_s) \prod_{t=j}^{n}f^2_{x_{t-1}}(x_t)\}\\
\label{theta1=theta2gen} f^{\theta_1=\theta_2>n}_x(\vec{x}_{1,n})&=& \rho\bar{\pi}p_1^n\prod_{s=1}^nf^0_{x_{s-1}}(x_s)\\
\label{nlesstheta1lesstheta2} f^{n<\theta_1<\theta_2}_x(\vec{x}_{1,n})&=& \bar{\rho}\bar{\pi}p_1^n\prod_{s=1}^nf^0_{x_{s-1}}(x_s).
\end{eqnarray}
Let us define the sequence of functions $S_n:\times_{i=1}^n\BbbE\rightarrow\Re$ as follows: $S_0(x_0)=1$ and for $n\geq 1$
\small
\begin{eqnarray}\label{disXuncond}
S_n(\vecx_n)
&=&f^{\theta_1\leq\theta_2\leq n}_x(\vec{x}_{1,n})+f^{\theta_1\leq n<\theta_2}_x(\vec{x}_{1,n})\\
\nonumber&&\mbox{}+f^{\theta_1=\theta_2>n}_x(\vec{x}_{1,n})+f^{n<\theta_1<\theta_2}_x(\vec{x}_{1,n}).
\end{eqnarray}

\begin{lemma}\label{recSn}
For $n>0$ the function $S_n(\vec{x}_{1,n})$ follows recursion
\begin{align}
\label{recurSn}S_{n+1}(\vec{x}_{1,n+1})&=\bH(x_n,x_{n+1},\vecPi_n)S_n(\vec{x}_{1,n})\\
\intertext{where\hfil}
\label{Hfunction}\bH(x,y,\alpha,\beta,\gamma)&= (1-\alpha) p_1 f^{0}_x(y)+ [p_2(\alpha -\beta) + q_1(1- \alpha-\gamma)]f^{1}_x(y)\\
\nonumber &  + [q_2\alpha+p_2\beta + q_1\gamma]f^{2}_x(y).
\end{align}
\end{lemma}
\begin{pf}
Let us assume $0\leq \theta_1\leq \theta_2$ and suppose that $B_i \in {\cB}$, $1\leq i \leq n+1$ and let us assume that $X_0=x$ and denote $D_n=\{\omega:X_i(\omega)\in B_i, 1\leq i\leq n\}$.
For $A_i=\{\omega:X_i\in B_i\}\in\cF_i$, $1\leq i\leq n+1$ we have by properties of the density function $S_n(\vecx_)$ with respect to the measure $\mu(\cdot)$ 
\begin{eqnarray*}
\int_{D_{n+1}}d\bP_x&=&\int_{\times_{i=1}^{n+1}B_i}S_{n+1}(\vecx_{n+1})\mu(d\vecx_{1,n+1})\\
&=&\int_{\times_{i=1}^nB_i}\int_{B_{n+1}}f(x_{n+1}|\vecx_{n})\mu(dx_{n+1})S_n(\vecx_{0,n})\mu(d\vecx_{1,n})\\
&=&\int_{\times_{i=1}^nB_i}\bP(A_{n+1}|\vecX_n=x_{n})\mu_x(d\vecx_{1,n})\\
&=&\int_{D_n}\bP_x(A_{n+1}|\vecX_{1,n})d\bP_x=\int_{D_n}\bP_x(A_{n+1}|\cF_n)d\bP_x=\int_{D_n}\one_{A_{n+1}}d\bP_x
\end{eqnarray*}
Now we split the conditional probability of $A_{n+1}$ into the following parts
\begin{eqnarray}\label{nProb1}
\bP_{x}(X_{n+1}\in A_{n+1}  \mid \cF_n ) &=& \bP_{x}(n<\theta_1<\theta_2,X_{n+1} \in A_{n+1} \mid \cF_n) \\
     &+& \bP_{x}(\theta_1 \leq n < \theta_2,X_{n+1} \in A_{n+1} \mid \cF_n) \label{nProb2}\\
     &+& \bP_{x}(n < \theta_1 =\theta_2,X_{n+1} \in A_{n+1} \mid \cF_n) \label{nProb3}\\
     &+& \bP_{x}(\theta_1\leq\theta_2 \leq n,X_{n+1} \in A_{n+1} \mid \cF_n)\label{nProb4}
\end{eqnarray}
\begin{description}
\item[In \eqref{nProb1}] we have:
\begin{eqnarray*}
\int_{D_n}\bP_{x}(\theta_2>\theta_1 > n,X_{n+1} \in A_{n+1} \mid \cF_n)d\bP_x &=&\int_{D_n}(\one_{\{\theta_1 = n+1\}} + \one_{\{\theta_1 >  n+1\}})\one_{A_{n+1}}d\bP_x\\
&\hspace{-26em}=&\hspace{-13em}\int_{\times_{i=1}^{n+1}B_i} (f^{n<\theta_1<\theta_2}_x(\vecx_{1,n}) (p_1f^0_{x_n}(x_{n+1})+q_1f^1_{x_n}(x_{n+1}))\mu(d\vecx_{1,n+1})\\
&\hspace{-26em}=&\hspace{-13em} \int_{\times_{i=1}^{n}B_i} (f^{n<\theta_1<\theta_2}_x(\vecx_{1,n}) \int_{B_{n+1}}(p_1f^0_{x_n}(x_{n+1})+q_1f^1_{x_n}(x_{n+1}))\mu(dx_{n+1}))\mu(d\vecx_{1,n})\\
&\hspace{-26em}=&\hspace{-13em}\int_{D_n}\bP_{x}(\theta_2>\theta_1 > n \mid \cF_n) [ \bP_{X_n}^0(A_{n+1})p_1+q_1\bP_{X_n}^{1}(A_{n+1}) ] d\bP_x.
\end{eqnarray*}

\item[In \eqref{nProb2}] we get by similar arguments as for (\ref{nProb1})
\begin{eqnarray*}
 \bP_{x}(\theta_1 \leq n < \theta_2 \!&,&\! X_{n+1} \in A_{n+1} \mid \cF_n) \nonumber\\
 &=& \bP_{x}(\theta_1 \leq n < \theta_2, \theta_2 = n+1, X_{n+1} \in A_{n+1} \mid \cF_n) \nonumber \\
 &&+ \bP_{x}(\theta_1 \leq n < \theta_2, \theta_2 \neq n+1, X_{n+1} \in A_{n+1} \mid \cF_n) \nonumber \\
 &=& \left( \bP_{x}(\theta_1 \leq n \mid \cF_n) - \bP_{x}(\theta_2 \leq n \mid \cF_n) \right)\nonumber\\
 && \times \! [ q_2 \bP_{X_n}^2(A_{n+1}) + p_2\bP_{X_n}^{1}(A_{n+1})] \nonumber
\end{eqnarray*}

\item[In \eqref{nProb4}] this part has the form:
\[
\bP_{x}(\theta_2 \leq n, X_{n+1} \in A_{n+1}\mid \cF_n) = \bP_{x}(\theta_2 \leq n \mid \cF_n)\bP_{X_n}^2(A_{n+1})
\]
\item[In \eqref{nProb3}] the conditional probability is equal to
\begin{eqnarray*}
 \bP_{x}(\theta_1=\theta_2>n \!&,&\! X_{n+1} \in A_{n+1} \mid \cF_n) \nonumber\\
 &=& \bP_{x}(\theta_1 = \theta_2>n, \theta_2 = n+1, X_{n+1} \in A_{n+1} \mid \cF_n) \nonumber \\
 &&+ \bP_{x}(\theta_1 = \theta_2>n, \theta_2  \neq n+1, X_{n+1} \in A_{n+1} \mid \cF_n) \nonumber \\
 &=& \bP_{x}(\theta_1 =\theta_2> n \mid \cF_n)[ q_1 \bP_{X_n}^2(A_{n+1}) + p_1\bP_{X_n}^{0}(A_{n+1})] \nonumber
\end{eqnarray*}
\end{description}
These formula lead to
\[
f(X_{n+1}|\vecX_{1,n})=\bH(X_n,X_{n+1},\Pi_n^1,\Pi_n^2,\Pi_n^{12}).
\]
which proves the lemma.
\end{pf}
\subsection{Conditional probability of various events defined by disorder moments}
According to definition of $\Pi^1_n$, $\Pi^2_n$, $\Pi^{12}_{n}$ we get
\begin{lemma}
\label{multiconddistform}
For the model discribed in Section~\ref{sformProblem} the following formulae are valid:
\begin{enumerate}
\item $\bP_x(\theta_2> \theta_1 > n|\cF_n)=1-\Pi_n^1-\Pi_n^{12} =\frac{f^{n<\theta_1<\theta_2}_x(\vec{x}_{1,n})}{S_n(\vecx_n)}$;
\item $\bP_x(\theta_2= \theta_1 > n|\cF_n)=\Pi_n^{12}=\frac{f^{\theta_1=\theta_2>n}_x(\vec{x}_{1,n})}{S_n(\vecx_n)}$;
\item $\bP_x(\theta_1\leq n<\theta_2|\cF_n)=\Pi^{1}_n-\Pi^2_n$;
\item $\bP_x(\theta_2\geq\theta_1>n|\cF_n)=1-\Pi_n^1 =\frac{\bar{\pi}p_1^n\prod_{s=1}^nf^0_{x_{s-1}}(x_s)}{S_n(\vecx_n)}$.
\end{enumerate}
\end{lemma}
\begin{pf}
\begin{enumerate}
\item We have 
\begin{eqnarray}\label{OMEGA1}
\Omega &=&\{\omega:n<\theta_1<\theta_2\}\cup\{\omega:\theta_1\leq n<\theta_2\}\\
\nonumber&&\cup\; \{\omega:\theta_1\leq\theta_2\leq n\}\cup\{\omega:\theta_1=\theta_2>n\}.
\end{eqnarray}
Hence $1=\bP_x(\omega:n<\theta_1<\theta_2|\cF_n)+(\Pi^1_n-\Pi^2_n)+\Pi^2_n+\Pi_n^{12}$ and
\begin{eqnarray*} 
\bP_x(\omega:n<\theta_1<\theta_2|\cF_n)&=&1-\Pi_n^1-\Pi_n^{12}.
\end{eqnarray*}

Let $B_i \in {\cB}$, $1\leq i \leq n$, $X_0=x$ and denote $D_n=\{\omega:X_i(\omega)\in B_i, 1\leq i\leq n\}$.
For $A_i=\{\omega:X_i\in B_i\}\in\cF_i$, $1\leq i\leq n$ and $D_{n}\in\cF_{n}$ we have 
\begin{eqnarray*}
\int_{D_n}\one_{\{\theta_2>\theta_1 > n\}}d\bP_x&=&\int_{D_n}\bP_x(\theta_2>\theta_1 > n|\cF_n)d\bP_x = \int_{D_n}\bP_x(\theta_2>\theta_1 > n|\vec{X}_n)d\bP_x\\
&=&\bP_x(\theta_2>\theta_1>n,D_n)=\int_{\times_{i=1}^nB_i} f^{n<\theta_1<\theta_2}_x(\vec{x}_{1,n})\mu(d\vecx_{1,n})\\
&=&\int_{\times_{i=1}^nB_i} f^{n<\theta_1<\theta_2}_x(\vec{x}_{1,n})(S_n(\vecx_n))^{-1}\mu_x(d\vecx_{1,n})\\
&=&\int_{D_n} f^{n<\theta_1<\theta_2}_x(\vec{X}_{1,n})(S_n(\vecX_n))^{-1}d\bP_x.
\end{eqnarray*}
Thus $\bP_x(\theta_2>\theta_1 > n|\cF_n) = \bar{\rho} \bar{\pi} p_1^n \prod_{i=1}^nf^0_{X_{i-1}}(X_i)(S_n(\vec{X}_n))^{-1}$.
\item  The second formula can be obtained by similar argument.
\item Let $\theta_1\leq\theta_2$. Since $\{\omega:\theta_2\leq n\}\subset\{\omega:\theta_1\leq n\}$ it follows that $\bP_x(\{\omega:\theta_1\leq n< \theta_n\}|\cF_n)=\bP_x(\{\omega:\theta_1\leq n\}\setminus\{\omega:\theta_2\leq n\}|\cF_n)=\Pi^1_n-\Pi^2_n$.
\end{enumerate}
These end the proof of the lemma.
\end{pf}

\begin{remark}\label{reccondprobab}
Let $B_i \in {\cB}$, $1\leq i \leq n+1$, $X_0=x$ and denote $D_n=\{\omega:X_i(\omega)\in B_i, 1\leq i\leq n\}$.
For $A_i=\{\omega:X_i\in B_i\}\in\cF_i$, $1\leq i\leq n$ and $D_{n}\in\cF_{n}$ we have 
\begin{eqnarray*}
\int_{D_n}\one_{\{\theta_1 > n\}}d\bP_x&=&\int_{D_n}\bP_x(\theta_1 > n|\cF_n)d\bP_x=\int_{D_n}\bP_x(\theta_1 > n|\vec{X}_n)d\bP_x\\
&=&\bP_x(\theta_1>n,D_n)=\int_{\times_{i=1}^nB_i} p_1^n\prod_{i=1}^nf^0_{x_{i-1}}(x_i)\mu(d\vecx_{1,n})\\
&=&\int_{\times_{i=1}^nB_i} p_1^n\prod_{i=1}^nf^0_{x_{i-1}}(x_i)(S_n(\vecx_n))^{-1}\mu_x(d\vecx_{1,n}).
\end{eqnarray*}
Thus $\bP_x(\theta_1 > n|\cF_n)=p_1^n\prod_{i=1}^nf^0_{X_{i-1}}(X_i)(S_n(\vec{X}_n))^{-1}$.
Moreover
\begin{eqnarray*}
1-\Pi^1_{n+1}&=&p_1f^0_{X_n}(X_{n+1})(1-\Pi^1_{n})S_{n}(\vec{X}_{n})(S_{n+1}(\vec{X}_{n+1}))^{-1}
\end{eqnarray*}
and $S_{n+1}(\vec{X}_{n+1})=\bH(X_n,X_{n+1},\vecPi^1_n)S_n(\vecX_n)$.
Hence
\[
\Pi^1_{n+1}=1-
\frac{p_1 f^0_{X_n}(X_{n+1})(1-\Pi^1_{n})}{\bH(X_n,X_{n+1},\vecPi_n)}.
\]
\end{remark}

\subsection{Some recursive formulae}
In derivation of the formulae in Theorem~\ref{reqform} the form of the distribution of some random vectors is taken into account. 
\begin{lemma}
\label{multidistform}
For the model discribed in Section~\ref{sformProblem} the following formulae are valid:
\begin{enumerate}
\item\label{loc1} $\bP_x(\theta_2=\theta_1>n+1|\cF_n)=p_1\Pi^{12}_n=p_1\rho(1-\Pi_n^1)$;
\item $\bP_x(\theta_2>\theta_1>n+1|\cF_n)=p_1(1-\Pi_n^1-\Pi_n^{12})$;
\item $\bP_x(\theta_1\leq n+1|\cF_n)=\bP_x(\theta_1\leq n+1<\theta_2|\cF_n)+\bP_x(\theta_2\leq n+1|\cF_n)$;
\item $\bP_x(\theta_1\leq n+1<\theta_2|\cF_n)=q_1(1-\Pi^{1}_n-\Pi^{12}_n)+p_2(\Pi^{1}_n-\Pi^{2}_n)$;
\item $\bP_x(\theta_2\leq n+1|\cF_n)=q_2\Pi^1_n+p_2\Pi_n^2+q_1\Pi^{12}_n$.
\item $\bP_x(\theta_1=m,\theta_2> n+1|\cF_n)=p_2\Pi_{m\;n}$.
\end{enumerate}
\end{lemma}
\begin{pf} 
\begin{enumerate}
\item On the set $D=\{\omega:X_0=x,X_1\in A_1,X_2\in A_2,\ldots,X_n\in A_n\}\in\cF_n$ we have
\begin{eqnarray*}
\int_D\one_{\{\theta_2=\theta_1>n+1\}}d\bP_x&=&\bP_x(D)\bP_x(\theta_2=\theta_1>n+1|D)\\
&=&\rho\bar{\pi}\sum_{j=n+2}^\infty p_1^{j-1}q_1\int_{\times_{i=1}^nA_i}\prod_{i=1}^n f^{0}_{x_{i-1}}(x_i)\mu(d\vecx_{1,n})\\
&=&p_1\rho\bar{\pi} p_1^{n}\int_{\times_{i=1}^nA_i}\prod_{i=1}^n f^{0}_{x_{i-1}}(x_i)\mu(d\vecx_{1,n})\\
&=&p_1\bP_x(D)\bP_x(\theta_2=\theta_1>n|D)=p_1\int_D\one_{\{\theta_2=\theta_1>n\}}d\bP_x.
\end{eqnarray*}
By (\ref{pi12nmx}) and the definition of the conditional probability this implies $\bP_x(\theta_2=\theta_1>n+1|\cF_n)=p_1\Pi^{12}_n$. Next,
\begin{eqnarray*}
\int_D\one_{\{\theta_1>n\}}d\bP_x&=&\bP_x(D)\bP_x(\theta_1>n|D)\\
&=&\bar{\pi}\sum_{j=n+1}^\infty p_1^{j-1}q_1\int_{\times_{i=1}^nA_i}\prod_{i=1}^n f^{0}_{x_{i-1}}(x_i)\mu(d\vecx_{1,n})\\
&=&\frac{1}{\rho}\bP_x(D)\bP_x(\theta_2=\theta_1>n|D)=\frac{1}{\rho}\int_D\one_{\{\theta_2=\theta_1>n\}}d\bP_x.
\end{eqnarray*}
These prove the part~1 of the lemma.
\item Similarly as above we get
\begin{eqnarray*}
\int_D\one_{\{\theta_2>\theta_1>n+1\}}d\bP_x&=&\bP(D)\bP_x(\theta_2>\theta_1>n+1|D)\\
&=&p_1\rho\bar{\pi} p_1^{n}\int_{\times_{i=1}^n A_i}\prod_{i=1}^n f^{0}_{x_{i-1}}(x_i)\mu(d\vecx_{1,n})\\
&=&p_1\bP(D)\bP_x(\theta_2>\theta_1>n|D)=p_1\int_D\one_{\{\theta_2>\theta_1>n\}}d\bP_x
\end{eqnarray*} 
By point~2 of Lemma~\ref{multiconddistform} we get the formula~2 of the lemma.
\item  It is obvious by assumption $\theta_1\leq\theta_2$.
\item  On the set $D$ we have
\begin{eqnarray*}
\int_D\one_{\{\theta_1\leq n+1<\theta_2\}}d\bP_x&=&\bP(D)\bP_x(\theta_1\leq n+1<\theta_2|D)\\
&\hspace{-16em}\stackrel{\eqref{rozkladyTeta},\eqref{rokladWarTeta2}}{=}&\hspace{-8em}\sum_{j=0}^{n+1}\bP(\omega:\theta_1=j)\sum_{k=n+2}^\infty\bar{\rho}p_2^{k-j-1}q_2\int_{\times_{i=1}^nA_i}\prod_{s=1}^{j-1}f^{0}_{x_{s-1}}(x_s)\prod_{r=j}^nf^{1}_{x_{r-1}}(x_r)\mu(d\vecx_{1,n})\\ 
&\hspace{-16em}=&\hspace{-8em}\bar{\pi}p_1^nq_1(1-\rho)\int_{\times_{i=1}^nA_i}\prod_{s=1}^{n}f^{0}_{x_{s-1}}(x_s)\mu(d\vecx_{1,n})\\
&\hspace{-16em}&\hspace{-8em}\mbox{\;}+p_2\sum_{0}^n\bP(\omega:\theta_1=j)p_2^{n+1-j}\int_{\times_{i=1}^nA_i}\prod_{s=1}^{j-1}f^{0}_{x_{s-1}}(x_s)\prod_{r=j}^nf^{1}_{x_{r-1}}(x_r)\mu(d\vecx_{1,n})\\
&\hspace{-16em}\stackrel{(L.\ref{multiconddistform})}{=}&\hspace{-8em}q_1\bP(D)\bP_x(\theta_2>\theta_1>n|D)+p_2\bP(D)\bP_x(\theta_1\leq n<\theta_2|D)\\
&\hspace{-16em}=&\hspace{-8em}q_1\int_D\one_{\{\theta_2>\theta_1>n\}}d\bP_x+p_2\int_D\one_{\{\theta_1\leq n<\theta_2\}}d\bP_x.
\end{eqnarray*}

\item If we substitute $n$ by $n+1$ in (\ref{OMEGA1}) than we obtain
\begin{eqnarray*}
\bP_x(\theta_2\leq n+1|\cF_n)&=&1-\bP_x(n+1<\theta_1=\theta_2|\cF_n)\\
&&-\bP_x(n+1<\theta_1<\theta_2|\cF_n)-\bP_x(\theta_1\leq n+1<\theta_2|\cF_n)\\
&=&1-p_1\Pi^{12}_n-p_1(1-\Pi^1_n-\Pi_n^{12})-q_1(1-\Pi^1_n-\Pi^{12}_n)\\
&&+p_2(\Pi^2_n-\Pi^1_n)
=q_2\Pi^1_n+p_2\Pi_n^2+q_1\Pi^{12}_n.
\end{eqnarray*}

\item We have
\begin{eqnarray*}
\int_D\one_{\{\theta_1=m,\theta_2>n+1\}}d\bP_x&=&\bP_x(D)\bP_x(\theta_1=m,\theta_2>n+1|D)\\
&\hspace{-16em}=&\hspace{-8em}\bar{\pi}\bar{\rho}p_1^{m-1}q_1\sum_{j=n+2}^\infty p_2^{j-m-1}q_2\int_{\times_{i=1}^nB_i}\prod_{i=1}^m f^{0}_{x_{i-1}}(x_i)\prod_{j=m+1}^n f^{1}_{x_{j-1}}(x_j)\mu(d\vecx_{1,n})\\
&\hspace{-16em}=&\hspace{-8em}p_2\bar{\pi}\bar{\rho} p_1^{m-1}q_1p_2^{n-m}\int_{\times_{i=1}^nB_i}\prod_{i=1}^m f^{0}_{x_{i-1}}(x_i)\prod_{j=m+1}^n f^{1}_{x_{j-1}}(x_j)\mu(d\vecx_{1,n})\\
&\hspace{-16em}=&\hspace{-8em}p_2\bP_x(D)\bP_x(\theta_1=m,\theta_2>n|D)=p_2\int_D\one_{\{\theta_1=m,\theta_2>n\}}d\bP_x.
\end{eqnarray*}
By (\ref{pinmx}) and the definition of conditional probability this implies $\bP_x(\theta_2=m,\theta_1>n+1|\cF_n)=p_2\Pi_{n\;m}$. 
These prove the part~6 of the lemma.
\end{enumerate}
\end{pf}

{\bf Acknowledgements.}{ I have benefited from discussions with Wojciech Sarnowski and Anna Karpowicz, for which I am grateful. I should like to thank Professor El\.zbieta Ferenstein for helpful suggestions. They provided numerous corrections to the manuscript.} 


\end{document}